\documentclass[11pt]{article}

\usepackage{graphicx}
\usepackage{subfigure}
\usepackage{color}
\usepackage{epsfig}
\usepackage{amssymb}
\usepackage{latexsym}
\usepackage{relsize}
\usepackage{amsmath,url}

\newcommand{\var}[1]{\mbox{$var [#1]$}}

\allowdisplaybreaks

\newcommand {\C} {{\rm I\kern-5.5pt C}}

\newcommand{\bP}[1]{{\mathbb{P}}\left[{#1}\right]}
\newcommand{\bE}[1]{{\mathbb{E}}\left[{#1}\right]}

\newcommand{\1}[1]{{\bf 1}\left[#1\right]}       

\newcommand{\fsquare}{\vrule height6pt width7pt depth1pt}   
\newcommand{\myproof}{{\hfill \\ \bf Proof. \ }}           
\newcommand{\myendpf}{\hfill\fsquare \\[0.1in]}             

\newtheorem{theorem}{Theorem}[section]

\newtheorem{lemma}[theorem]{Lemma}
\newtheorem{proposition}[theorem]{Proposition}
\newtheorem{corollary}[theorem]{Corollary}

\begin{document}

\title{Zero-one laws for connectivity in inhomogeneous random key graphs}

\author{
Osman Ya\u{g}an \\
{\tt oyagan@ece.cmu.edu}\\
Department of Electrical and Computer Engineering and CyLab\\
Carnegie Mellon University\\
}

\maketitle

\begin{abstract}
\normalsize 
We introduce a new random key predistribution scheme
for securing {\em heterogeneous} wireless sensor networks. Each of the $n$
sensors in the network is classified into 
$r$ classes according to some probability distribution $\boldsymbol{\mu}=\{\mu_1,\ldots,\mu_r\}$. 
Before deployment, a class-$i$ sensor is assigned $K_i$ cryptographic keys that are selected uniformly at random from a common pool of $P$ keys. Once deployed, a pair of sensors can communicate
{\em securely} if and only if they have a key in common. We model the communication topology of this network by a newly defined {\em inhomogeneous} random key graph. We establish scaling conditions on the parameters $P$ and $\{K_1,\ldots, K_r\}$ so that this graph i) has no isolated nodes; and ii) is connected, both with high probability. The results are given in the form of zero-one laws with the number of sensors $n$ growing unboundedly large; {\em critical} scalings are identified and shown to coincide for both graph properties. Our results 
are shown to complement and improve those given by Godehardt et al. and Zhao et al. for the same model, therein referred to as the {\em general random intersection graph}.
 
\end{abstract}


\section{Introduction}\label{introduction}

Random key graphs are naturally induced by 
the Eschenauer-Gligor (EG) random key predistribution scheme \cite{virgil},
which is a widely recognized solution
for securing
wireless sensor network (WSN) communications \cite{adrian,1354530}. 
Denoted by $\mathbb{G}(n, K, P)$, random key graph is constructed on the
vertices $\mathcal{V}=\{v_1, v_2, \ldots, v_n\}$ as follows. 
Each vertex $v_i$ is assigned {\em independently} a set $\Sigma_i$ 
of $K$ cryptographic keys that
are picked uniformly at random from a pool of size $P$. Then, any pair of 
vertices $v_i, v_j$ are
{\em adjacent} if they share a key, i.e.,
if $\Sigma_i \cap \Sigma_j \neq \emptyset$. 
Random key graphs have recently received attention in a wide range of areas including
modeling small world networks \cite{5383986}, recommender systems \cite{Marbach}, and 
clustering and classification analysis \cite{GodehardtJaworskiRybarczyk}.
Properties that have been studied 
include absence of isolated nodes \cite{YM_ISIT2008}, connectivity \cite{yagan,ryb3}, $k$-connectivity \cite{ZYG_IT_k}, and $k$-robustness \cite{ZhaoCDC},
among others.

In this paper we propose and study a variation of the EG
scheme that is more suitable for {\em heterogeneous} WSNs;
it is in fact envisioned that many military and commercial WSN applications will consist of
heterogeneous nodes \cite{het_2,het_1}.
Namely, we assume that the network consists of sensors with varying level of resources 
(e.g., computational, memory, power) and possibly with varying level of security and connectivity 
requirements. As a result of this heterogeneity, it may no longer be sensible to assign the same number
of keys to all sensors in the network as prescribed by the EG scheme. 
Instead, we consider a scheme
where the number of keys that will be assigned 
to each sensor is independently drawn from the set $\boldsymbol{K}=\{K_1,\ldots, K_r\}$
according to some probability distribution $\boldsymbol{\mu}=\{\mu_1,\ldots,\mu_r\}$,
for some fixed integer $r$. We can think of this as each vertex $v_x$ being
assigned to a priority class-$i$ with probability $\mu_i>0$ and then receiving 
a key ring with the size $K_i$ associated with this class. As before, we assume that 
once its size is fixed, the key ring $\Sigma_x$ is constructed by sampling the key pool
randomly and without replacement.

Let $\mathbb{G}(n;\boldsymbol{\mu},\boldsymbol{K},P)$ denote the 
random graph induced by the heterogeneous key predistribution 
scheme described above, where again 
a pair of nodes are adjacent  as long as they share a key; see Section \ref{sec:NetworkModel}
for precise definitions.
Inspired by the recently studied inhomogeneous Erd\H{o}s-R\'enyi (ER) graphs
\cite{BollobasJansonRiordan,devroye2014connectivity},
we refer to this graph as the {\em inhomogeneous} random key graph.
The main goal of this paper is to study connectivity properties
of $\mathbb{G}(n;\boldsymbol{\mu},\boldsymbol{K},P)$ and
to understand how the parameters $n, \boldsymbol{\mu},\boldsymbol{K},P$ should behave
so that the resulting graph is connected almost surely. Such results can
be useful in deriving guidelines for designing
heterogenous WSNs so that they are securely connected. By comparison with 
the results for the standard random key graph, they can also shed light
on the effect of heterogeneity on the connectivity properties of WSNs.

We establish zero-one laws for the property that $\mathbb{G}(n;\boldsymbol{\mu},\boldsymbol{K},P)$
has no {\em isolated} nodes (see Theorem \ref{thm:node_isolation}) and for
the property that  $\mathbb{G}(n;\boldsymbol{\mu},\boldsymbol{K},P)$ is connected (see Theorem \ref{thm:connectivity}). 
Namely, we scale the parameters $\boldsymbol{K}$ and $P$ and provide critical conditions
on this scaling such that the resulting graph almost surely has no isolated node (resp.~connected) and almost surely has at least one isolated node (resp.~connected), respectively, when the number of nodes $n$ goes to infinity. 
We show that the critical conditions for the two graph properties coincide, meaning that 
absence of isolated nodes and connectivity are asymptotically equivalent properties for the inhomogeneous random
key graph.
Other well-known models that exhibit the same behavior include  ER graphs \cite{bollobas}, random key graphs \cite{yagan}, and random geometric graphs \cite{penrose_book}. 

Our results are also compared with the existing results by Zhao et al. \cite{ZhaoCDC} and Godehardt et al. \cite{GodehardtJaworskiRybarczyk} for the $k$-connectivity and connectivity, respectively, of $\mathbb{G}(n;\boldsymbol{\mu},\boldsymbol{K},P)$; in those references $\mathbb{G}(n;\boldsymbol{\mu},\boldsymbol{K},P)$
was referred to as a {\em general} random intersection graph. We show that earlier results are constrained
to parameter ranges 
that are unlikely to be feasible in real world WSN implementations 
due to excessive memory requirement or very limited resiliency
against adversarial attacks. On the contrary, our results cover parameter ranges that are widely regarded as feasible
for most WSNs; see Section \ref{sec:comparison} for details.

In addition, our main results indicate that the minimum key ring size in the network has a significant impact on the connectivity of $\mathbb{G}(n;\boldsymbol{\mu},\boldsymbol{K},P)$, perhaps in a way that would be deemed surprising. In particular, for the standard random key graph $\mathbb{G}(n;K,P)$ the critical threshold 
for connectivity and absence of isolated nodes is known \cite{yagan,ryb3} to be given by  $\frac{K^2}{P}\sim c \frac{\log n}{n}$ and the resulting graph is asymptotically almost surely connected (resp.~not connected) if $c>1$ (resp.~$c<1$). 
For the inhomogeneous random key graph $\mathbb{G}(n;\boldsymbol{\mu},\boldsymbol{K},P)$ one would be tempted to think that an equivalent result holds under the scaling $\frac{K_{\textrm{avg}}^2}{P}\sim  c \frac{\log n}{n}$, with $K_{\textrm{avg}}=\sum_{j=1}^{r}\mu_j K_j$ denoting the mean key ring size. Instead, we show that
the zero-one laws for absence of isolated nodes and connectivity hold under
$
\frac{K_{\textrm{min}}K_{\textrm{avg}}}{P}\sim c \frac{\log n}{n}$,
where $K_{\textrm{min}}$ stands for the minimum of $\{K_1,\ldots, K_r\}$; see Corollary \ref{cor:node_isolation}.
This implies that in the heterogeneous key predistribution scheme, the mean number of keys required per sensor node to achieve connectivity can be significantly larger than that required in the homogeneous case. For instance, the expense of allowing an arbitrarily small fraction of sensors to keep half as many keys as in the homogeneous case
would be to increase the average key ring size by two-fold.

The rest of the paper is organized as follows. In
Section \ref{sec:NetworkModel}, we give detailed description of the heterogeneous 
random key predistribution scheme and the resulting inhomogeneous random key graph model.
Section \ref{sec:main:res}
is devoted to presenting the main results of the paper, namely the zero-one
laws for absence of isolated nodes (see Theorem \ref{thm:node_isolation}) and
for connectivity (see Theorem \ref{thm:connectivity}) in inhomogeneous random key graphs. There,
we also compare our results with relevant work from the literature and also comment
on the implications of our results on designing secure sensor networks.
In Section \ref{sec:prelim}, we present some preliminary technical results  
that will be useful in proving the main results of the paper. 
The proof of Theorem \ref{thm:node_isolation} is outlined in Section 
\ref{sec:proof_thm_1} with necessary technical steps completed in
Section \ref{sec:prop_proofs}. We start the proof of Theorem
\ref{thm:connectivity} in Section \ref{sec:proof_thm_2}
and complete it in Sections \ref{sec:OneLawAfterReductionPart1} through 
Section \ref{sec:OneLawAfterReductionPart3}. 

We close with a word on notation and conventions in use. All limiting
statements, including asymptotic equivalences, are understood with
the number of sensor nodes $n$ going to infinity. The random
variables (rvs) under consideration are all defined on the same
probability triple $(\Omega, {\cal F}, \mathbb{P})$. Probabilistic
statements are made with respect to this probability measure
$\mathbb{P}$, and we denote the corresponding expectation and variance
operators
by $\mathbb{E}$ and {\mbox{$var$}}, respectively.
We use the notation $=_{st}$ to indicate
distributional equality. 
The indicator function of an event $E$ is
denoted by $\1{E}$, while $E^c$ denotes the complement
of $E$.
We say that an even holds {\em with high
probability} (whp) if it holds with probability $1$ as $n \to
\infty$. 
For any discrete set $S$ we write $|S|$ for its
cardinality.
 In comparing
the asymptotic behaviors of the sequences $\{a_n\},\{b_n\}$,
we use
$a_n = o(b_n)$,  $a_n=w(b_n)$, $a_n = O(b_n)$, $a_n = \Omega(b_n)$, and
$a_n = \Theta(b_n)$, with their meaning in
the standard Landau notation. We also use $a_n \sim b_n$ to denote the
asymptotic equivalence $\lim_{n \to \infty} {a_n}/{b_n}=1$.


\section{Model Definitions}
\label{sec:NetworkModel}

Consider a network that consists of $n$ sensor nodes labeled as $v_1, \ldots, v_n$.
Our key predistribution idea is based on classifying the nodes in this network into $r$ sets 
(e.g., depending on their level of importance) and then assigning different number of cryptographic keys to sensors 
based on their class.
Assume that each of the $n$ nodes in the network are independently assigned to a class according to some probability distribution $\boldsymbol{\mu}:\{1,\ldots, r\} \to (0,1)$. Namely, with $t_x$ denoting the class (or, type) of node $v_x$, we have
\[
\bP{{t_{\ell}=i}} = \mu_{i}>0, \qquad 
  i = 1,\ldots, r,
\] 
for each $\ell=1,\ldots, n$.
Then, a class-$i$ node is assigned $K_{i}$ keys that are selected uniformly at random from a pool of size $P$, for each $i=1,\ldots, r$.  More precisely, the key ring $\Sigma_x$ of a node $x$ is an $\mathcal{P}_{K_{t_x}} $-valued random variable (rv) where
$\mathcal{P}_{A} $ denotes the collection of all subsets of $\{ 1,
\ldots , P \}$ which contain exactly $A$ elements -- Obviously, we
have $|\mathcal{P}_{A} | =  {P \choose A}$. It is further assumed that the rvs $\Sigma_1,
\ldots , \Sigma_n$ are independent and identically distributed. 

Let $\boldsymbol{K}=(K_1, \ldots, K_r)$ and $\boldsymbol{\mu}=(\mu_1, \ldots, \mu_r)$. Without loss of generality we assume that $K_1 \leq K_2 \leq \cdots \leq K_r$. Consider a random graph $\mathbb{G}$ defined on the vertex set $\mathcal{V}=\{v_1,\ldots,v_n\}$ such that 
two nodes $v_x$ and $v_y$ are adjacent, denoted $v_x \sim v_y$, if they have at least one key in common in their corresponding key rings.  Namely, we have
\begin{equation}
v_x \sim v_y \qquad \textrm{if} \qquad \Sigma_x \cap \Sigma_y \neq \emptyset.
\label{eq:adjacency}
\end{equation}

The adjacency condition (\ref{eq:adjacency})
defines the inhomogeneous random key graph, hereafter denoted $\mathbb{G}(n;\boldsymbol{\mu}, \boldsymbol{K},P)$. The name is reminiscent of the recently studied inhomogeneous random
graph \cite{BollobasJansonRiordan} model where nodes are again divided into $r$ classes, and a
class $i$ node and a class $j$ node are connected with probability $p_{ij}$ independent of everything else. This independence disappears in the inhomogeneous random key graph case, but one can still compute $p_{ij}$ as
\begin{equation}
p_{i j} := 1 -{{{P-K_i}\choose{K_j}}\over{{P}\choose{K_j}}}, \quad i,j=1,\ldots, r.
\label{eq:edge_prob}
\end{equation}
In view of (\ref{eq:edge_prob}), our key predistribution scheme results in higher priority nodes (i.e., nodes with more assigned keys) connecting with each other with higher probability; see Proposition \ref{prop:ordering_lamda}.
In presenting our results below, we shall make use of the {\em mean} probability of edge occurrence for each node class. Namely, we define
\begin{equation}
\lambda_i := \sum_{j=1}^{r} p_{ij}  \mu_j, \quad i=1, \ldots, r,
\label{eq:mean_edge_prob}
\end{equation}
where $p_{ij}$ denotes the probability that a node of class-$i$ and a node of class-$j$ have an edge in between; see (\ref{eq:edge_prob}). It is easy to see that the mean number of edges incident on a node (i.e.,  the {\em degree} of a node) of class-$i$ is given by $(n-1)\lambda_i$.

Throughout, we assume that the number of classes $r$ is fixed and do not scale with $n$, and so are the probabilities $\mu_1, \ldots, \mu_r > 0$. All remaining parameters are assumed to be scaled with $n$, and we shall be interested in the properties of the resulting inhomogeneous random key graph as $n$ grows unboundedly large. The dependence
of scheme parameters and events on $n$ will be denoted by a subscript, while 
that of some variables will be denoted by a parenthesis. For instance, we define
\begin{equation}
p_{i j} (n) := 1 -{{{P_n-K_{i,n}}\choose{K_{j,n}}}\over{{P_n}\choose{K_{j,n}}}}, \quad i,j=1,\ldots, r,
\label{eq:edge_prob_n}
\end{equation}	
and 
\begin{equation}
\lambda_i(n) := \sum_{j=1}^{r} p_{ij}(n)  \mu_j, \quad i=1, \ldots, r.
\label{eq:mean_edge_prob_n}
\end{equation}


\section{Main Results and Discussion} \label{sec:main:res}
\subsection{The results}
Our main results are presented next. To fix the terminology,
we refer to any mapping $K_1,\ldots,K_r, P: \mathbb{N}_0 \to \mathbb{N}_0^{r+1}$ as a {\em scaling} as long as the 
conditions
\begin{equation}
1\leq K_{1,n} \leq K_{2,n} \leq \cdots \leq K_{r,n} < P_n
\label{eq:scaling_def}
\end{equation}
are satisfied for all $n=2,3, \ldots$. To simplify the notation, we also let $\boldsymbol{K}_n=(K_{1,n}, K_{2,n}, \ldots,  K_{r,n})$.
We first present a zero-one law for the absence of isolated nodes in inhomogeneous random key graphs.
\begin{theorem}
{\sl Consider a probability distribution $\boldsymbol{\mu}=(\mu_1,\ldots,\mu_r)$ with $\mu_i>0$ for
 $i=1,\ldots, r$, and a scaling  $K_1,\ldots,K_r, P: \mathbb{N}_0 \to \mathbb{N}_0^{r+1}$ such that 
\begin{equation}
\lambda_1(n) \sim c \frac{\log n}{n}
\label{eq:scaling_law}
\end{equation}
for some $c>0$. Then, we have
\begin{align}
 \lim_{n \to \infty}\bP{\mathbb{G}(n;\boldsymbol{\mu}, \boldsymbol{K}_n,P_n) ~\textrm{has no isolated nodes}~} 
&= \left\{
\begin{array}{cc}
0  & \textrm{if ~ $c<1$}    \\
   &      \\
1 &      \textrm{if ~ $c>1$.}   
\end{array}
\right.
\label{eq:zero_one_law}
\end{align}
}
\label{thm:node_isolation}
\end{theorem}
 
A proof of Theorem \ref{thm:node_isolation} can be found in Section \ref{sec:proof_thm_1}.
The scaling condition (\ref{eq:scaling_law}) will often be used in the equivalent form 
\begin{equation}
\lambda_1(n) = c_n \frac{\log n}{n}
\label{eq:scaling_law_2}
\end{equation}
with $\lim_{n \to \infty} c_n = c > 0$.

Next, we present an analogous result for the property of graph connectivity.

\begin{theorem}
{\sl Consider a probability distribution $\boldsymbol{\mu}=(\mu_1,\ldots,\mu_r)$ with $\mu_i>0$ for
 $i=1,\ldots, r$, and a scaling  $K_1,\ldots,K_r, P: \mathbb{N}_0 \to \mathbb{N}_0^{r+1}$ such that 
(\ref{eq:scaling_law}) holds
for some $c>0$. Under the assumptions
\begin{align}
P_n = \Omega(n)
\label{eq:condition_for_con_1}
\end{align}
and
\begin{align}
\frac{(K_{1,n})^2}{P_n} = w\left(\frac{1}{n}\right),
\label{eq:condition_for_con_2}
\end{align}
we have
\begin{align}
 \lim_{n \to \infty}\bP{\mathbb{G}(n;\boldsymbol{\mu}, \boldsymbol{K}_n,P_n) ~\textrm{is connected}~} 
&= \left\{
\begin{array}{cc}
0  & \textrm{if ~ $c<1$}    \\
   &      \\
1 &      \textrm{if ~ $c>1$.}   
\end{array}
\right.
\label{eq:zero_one_law_con}
\end{align}
}
\label{thm:connectivity}
\end{theorem}

The condition (\ref{eq:condition_for_con_1}) implies that there exists a constant $\sigma>0$ such 
that 
\begin{align}
P_n &\geq \sigma n
\label{eq:condition_for_con_1b}
\end{align}
for all $n=2,3,\ldots$  sufficiently large.

In words, Theorem \ref{thm:node_isolation} (resp.~Theorem \ref{thm:connectivity})
states that the inhomogeneous random key 
graph $\mathbb{G}(n;\boldsymbol{\mu}, \boldsymbol{K}_n,P_n)$ 
has no isolated node (resp.~is connected) whp if the mean degree of \lq\lq the nodes that have the least number of keys" is scaled as $(1+\epsilon) \log n$ for some $\epsilon>0$; in view of Proposition \ref{prop:ordering_lamda}, the nodes that are assigned the least number of keys have the {\em minimum} mean-degree in the graph.
On the other hand, if this minimal mean degree scales like $(1-\epsilon) \log n$ for some $\epsilon>0$, then whp $\mathbb{G}(n;\boldsymbol{\mu}, \boldsymbol{K}_n,P_n)$ has a node that is isolated, and hence not connected. 
The additional conditions (\ref{eq:condition_for_con_1}) and (\ref{eq:condition_for_con_2})
are enforced here merely for technical reasons and
are required only for the one-law part of the connectivity result, Theorem \ref{thm:connectivity}. 
A detailed discussion on these additional conditions is given in Section \ref{sec:comments_conditions}, where 
we explain why they are likely to hold in many real-world WSN applications. There, we also discuss how 
and when these conditions can be relaxed or replaced by milder conditions.

Our results demonstrate that the inhomogeneous random key graph provides one more example random graph model where the properties of absence of isolated nodes and connectivity are asymptotically equivalent. Other well-known examples include Erd\H{o}s-R\'enyi graphs \cite{bollobas}, random key graphs \cite{yagan}, random geometric graphs \cite{penrose_book}, and intersection of random key graphs and ER graphs \cite{yagan_onoff}.
Our results are also analogous to the recent findings by Levroye and Freiman \cite{devroye2014connectivity}
for the connectivity of inhomogeneous Erd\H{o}s-R\'enyi graph model, where nodes are classified into $r$ classes independently according to a probability distribution $\boldsymbol{\mu}$ and an edge is drawn between a class-$i$ and a class-$j$ node with probability $p_{ij}(n)$ independent of everything else.
With $\lambda_i(n)$ defined as in (\ref{eq:mean_edge_prob_n}), their result states that if $\min_{i=1,\ldots,r}\lambda_i(n) \sim c \log n/n$ then with $c>1$ (resp.~$c<1$) the corresponding graph is connected (resp.~not connected) whp, under some additional technical conditions.\footnote{Results in \cite{devroye2014connectivity} cover more general cases than presented here; e.g., the case where the number of classes $r$ is not bounded.} 

We now present a corollary that states the zero-one laws of Theorem \ref{thm:node_isolation} and Theorem \ref{thm:connectivity} 
under a different scaling condition than (\ref{eq:scaling_law}). 
This alternative formulation will make it easier to derive design guidelines 
for  {\em dimensioning} heterogeneous key predistribution schemes,
namely in adjusting key ring sizes $K_1,\ldots,K_r$  and probabilities $\mu_1,\ldots,\mu_r$ such that
the resulting network i) has no isolated sensors and ii) is connected, both whp.

\begin{corollary}
{\sl Consider a probability distribution $\boldsymbol{\mu}=(\mu_1,\ldots,\mu_r)$ with $\mu_i>0$ for
 $i=1,\ldots, r$ and a scaling  $K_1,\ldots,K_r, P: \mathbb{N}_0 \to \mathbb{N}_0^{r+1}$. Let $|\Sigma|_n$
denote a rv that takes the value $K_{i,n}$ with probability $\mu_i$ for each $i=1,\ldots, r$. 
If it holds that
\begin{equation}
\frac{K_{1,n} \bE{|\Sigma|_n}}{P_n} \sim c \frac{\log n}{n}
\label{eq:scaling_law_cor}
\end{equation}
for some $c>0$, then we have the zero-one law
(\ref{eq:zero_one_law}) for absence of isolated nodes. If, in addition, the
the conditions (\ref{eq:condition_for_con_1}) and (\ref{eq:condition_for_con_2}) are satisfied by this scaling, then 
we also have the zero-one law (\ref{eq:zero_one_law_con}) for connectivity.
}
\label{cor:node_isolation}
\end{corollary}
A proof of Corollary \ref{cor:node_isolation} is given in Appendix \ref{sec:proof_cor}, where we
show that the scaling conditions (\ref{eq:scaling_law}) and (\ref{eq:scaling_law_cor}) are indeed {\em equivalent} to each other, meaning that 
one can obtain both Theorem \ref{thm:node_isolation} and Theorem \ref{thm:connectivity}
from Corollary \ref{cor:node_isolation}, and vice versa.
 We remark that
$\bE{|\Sigma|_n}$ gives the mean number of keys assigned to a sensor in the network. With this in mind,
Corollary \ref{cor:node_isolation} provides various design choices to ensure that no sensor is isolated 
in the network: One just has to set the minimum and average key ring sizes such that their multiplication 
scales as $(1+\epsilon) \frac{P_n \log n}{n}$ for some $\epsilon>0$. We also see from
Corollary \ref{cor:node_isolation} that such a scaling would also ensure connectivity as long as the
additional conditions (\ref{eq:condition_for_con_1})-(\ref{eq:condition_for_con_2}) are also satisfied. 

To compare with the homogeneous random key predistribution scheme,
set $r=1$ and consider a universal key ring size $K_n$ in Corollary \ref{cor:node_isolation}. 
This leads to zero-one laws for the absence of isolated nodes and connectivity in the
standard random key graph $\mathbb{G}(n;K_n,P_n)$. Namely, with 
\begin{equation}
\frac{K_n^2}{P_n} \sim c \frac{\log n}{n}, \qquad c>0
\label{eq:scaling_rkg}
\end{equation}
analogs of (\ref{eq:zero_one_law}) and (\ref{eq:zero_one_law_con}) 
are obtained for $\mathbb{G}(n;K_n,P_n)$; these results had
already been established \cite{YM_ISIT2008,yagan} by the authors (in stronger forms).
An interesting observation is that minimum key ring size has a dramatic impact on the 
connectivity properties of inhomogeneous random key graph. To provide a simple and concrete example,
set $P_n = n \log n$. In the homogeneous case, we see from (\ref{eq:scaling_rkg}) that the universal key ring size
has to scale as $K_n = (1+\epsilon) \log n$ for some $\epsilon>0$ to ensure that the network is free of isolated nodes
and is connected.
In the heterogeneous case, one gains the flexibility of having a positive fraction of sensors in the network with substantially smaller number of keys. For the absence of isolated nodes, a positive fraction of sensors can be assigned as few as one key per node. However, from Corollary \ref{cor:node_isolation} we see that this comes at the expense of having to assign a substantially larger key rings 
to a positive fraction of other sensors
in the network. More precisely,  if $K_{1,n}=O(1)$ then we must have $K_{r,n}=\Omega((\log n)^2)$
to have no isolated nodes under the same setting. For connectivity on the other hand, 
we see from (\ref{eq:milder_extra_condition}) that the minimum 
key ring size $K_{1,n}$ can be kept on the order  of  $O(\sqrt{\log n})$ and connectivity
can still be achieved with mean key ring size satisfying $O((\log n)^{1.5})$.


\subsection{Comments on the technical conditions (\ref{eq:condition_for_con_1})-(\ref{eq:condition_for_con_2}) of Theorem \ref{thm:connectivity}}
\label{sec:comments_conditions}

We now provide a detailed discussion on the technical conditions 
(\ref{eq:condition_for_con_1}) and (\ref{eq:condition_for_con_2}) enforced in Theorem \ref{thm:connectivity}.
We will focus on i) the feasibility of these additional conditions for real-world WSN implementations, and ii)
when and how they can be replaced with milder conditions. 

We start with the condition (\ref{eq:condition_for_con_1}) that states the key pool size grows at least linearly
with the network size $n$. In terms of applicability in the context of heterogeneous key predistribution schemes in WSNs, this condition is not stringent at all. In fact, it is often needed that key pool size $P_n$ be much larger than the network size $n$ \cite{DiPietroMeiManciniPanconesiRadhakrishnan2004,DiPietroTissec,virgil} as otherwise the network will be extremely vulnerable against node capture attacks. From a technical point of view, the case where $P_n = \Omega(n)$ is also the more interesting and challenging one as compared to the case where $P_n=o(n)$. For instance, when 
$P_n=O(n^{\delta})$ for some $0< \delta < 1/2$,
the inhomogeneous random key graph $\mathbb{G}(n;\boldsymbol{\mu}, \boldsymbol{K}_n,P_n)$ can be shown to be  connected for any $\boldsymbol{\mu}$ as long as $K_{1,n} \geq  2$ ; see \cite[Lemma 8.1]{YM_small_pool} for a proof of a similar result for the standard random key graph. This means that if $P_n=O(n^{\delta})$ with $\delta<1/2$, even two keys per sensor node is enough to get network connectivity whp. Finally, we remark that the scaling condition
(\ref{eq:scaling_law}) or its equivalent (\ref{eq:scaling_law_cor}) already implies that 
$
P_n = \Omega (\frac{n}{\log n})
$
since $K_{1,n} \bE{|\Sigma|_n}\geq 1$.

Next, we look at the condition (\ref{eq:condition_for_con_2}) and start with 
discussing possible relaxations. First of all, (\ref{eq:condition_for_con_2}) is stronger
than what is actually needed for our proof to work; it is enforced to enable a shorter proof and an easier
exposition of the main result.
By inspection of the arguments in Section \ref{subsec:range4}, it can be seen that (\ref{eq:condition_for_con_2}) can be replaced with
\begin{align}
\left\{
\begin{array}{ll}
  \frac{K_{1,n}^2}{P_n} \geq \frac{\frac{2 \log 2 + \log (1-\mu_r)+\epsilon}{\beta \nu}}{n} \quad \textrm{and} \quad K_{1,n}=w(1),   & \textrm{if $\mu_r \leq 0.75$}  \\
  & \\
\frac{K_{1,n}^2}{P_n} =\Omega \left( \frac{1}{n (\log n)^M} \right) \quad \textrm{and} \quad K_{1,n}=w(1),        &   \textrm{if $\mu_r > 0.75$}
\end{array}
\right.
\label{eq:milder_extra_condition}
\end{align}
with some $\beta>0$ and $\nu>0$ to be specified, and for any $\epsilon>0$ and any finite integer $M$; the details are omitted here for brevity. 

As we look at (\ref{eq:milder_extra_condition}), we see that $K_{1,n}=w(1)$  is needed for any $\mu_r$.
In fact, this condition can easily be satisfied in real-world WSN implementations given that
key ring sizes on order of $O(\log n)$ are regarded as feasible for most sensor networks \cite{DiPietroTissec}. 
Considered in combination with (\ref{eq:scaling_law_cor}), other conditions enforced in 
(\ref{eq:milder_extra_condition}) bounds the {\em variability} in the key ring sizes used in the network. In particular, given that
\[
\frac{\bE{|\Sigma|_n}}{K_{1,n}} = \frac{K_{1,n} \bE{|\Sigma|_n}}{P_n} \left( \frac{(K_{1,n})^2 }{P_n} \right)^{-1}
= \Theta \left(\frac{\log n}{n} \right) \left( \frac{(K_{1,n})^2 }{P_n} \right)^{-1},
\]
(\ref{eq:milder_extra_condition})  implies
$\frac{\bE{|\Sigma|_n}}{K_{1,n}} = O\left(\log n\right)$ when $\mu_r \leq 0.75$ and 
$\frac{\bE{|\Sigma|_n}}{K_{1,n}} = O\left((\log n)^{M}\right)$ when $\mu_r > 0.75$.
Thus, we see that when more than 75 \% of the sensors receive the largest key rings, one can afford to use much smaller key rings for the remaining sensors, as compared to the case when $\mu_r \leq 0.75$. 

Collecting, while conditions enforced in (\ref{eq:milder_extra_condition}) take away from the flexibility of assigning {\em very small} key rings to a certain fraction of sensors (as we were allowed to do for the absence of isolated nodes),
they can still be satisfied easily in most real-world implementations. To provide a concrete example, one can set $P_n=n \log n$ and have $K_{1,n}= (\log n)^{1/2+\epsilon}$ and $\bE{|\Sigma|_n}=(1+\epsilon)(\log n)^{3/2-\epsilon}$ with any $\epsilon > 0$; in view of Theorem \ref{thm:connectivity} and (\ref{eq:milder_extra_condition}) the resulting network will be connected whp. With the same
$P_n$, it is possible to have much smaller $K_{1,n}$ when $\mu_r>0.75$. For example, we can have
$K_{1,n}=\log \log \cdots \log n$ and $\bE{|\Sigma|_n}=\Omega((\log n)^2)$.
Of course, one can also have all key ring sizes on the same order and set $K_{1,n} = c_1 \log n$ and
$\bE{|\Sigma|_n} =c_2 \log n$ with $c_1 c_2 > 1$, to obtain a connected WSN whp.


 \subsection{Comparison with related work}
\label{sec:comparison}

The random graph model $\mathbb{G}(n;\boldsymbol{\mu}, \boldsymbol{K}_n,P_n)$ considered here is 
also known as {\em general random intersection graphs} in the literature; e.g., see \cite{ZhaoCDC,Rybarczyk,GodehardtJaworski}. To the best of our knowledge this model has been first considered by Godehardt and Jaworski \cite{GodehardtJaworski} and by Goderhardt et al. \cite{GodehardtJaworskiRybarczyk}. Results for both the existence of isolated nodes and graph connectivity have been established; see below for a comparison of these results with those established here. Later, Bloznelis et al. \cite{Rybarczyk} analyzed the component evolution problem in the general random intersection graph and provided scaling conditions for the existence of a {\em giant component}. There, they also established that under certain conditions
$\mathbb{G}(n;\boldsymbol{\mu}, \boldsymbol{K}_n,P_n)$ behaves very similarly with a standard Erd\H{o}s-R\'enyi graph \cite{bollobas}. Taking advantage of this similarity,  Zhao et al. \cite{ZhaoCDC} 
established various results for the $k$-connectivity and $k$-robustness of the general random intersection graph by means of a coupling argument.

We now compare our results with those established in the literature. Our main argument is that 
previous results for the connectivity of inhomogeneous random key graphs are constrained to very narrow parameter ranges that are impractical
for wireless sensor network applications. In particular, we will argue below that the result by Zhao et al. \cite{ZhaoCDC} is restricted to {\em very large} key ring sizes, rendering them impractical for resource-constrained sensor networks. On the other hand, the results by Godehardt et al. \cite{GodehardtJaworski, Rybarczyk} focus on fixed key ring sizes that do not grow with the network size $n$. As a consequence, in order to ensure connectivity, their result requires a key pool size $P_n$ that is {\em much smaller} than typically prescribed for security and resiliency purposes.

To fix the terminology, let $\mathcal{D}_n:\{1,2,\ldots,P_n\} \to [0,1]$ be the probability distribution used for drawing the {\em size} of the key rings $\Sigma_1,\ldots, \Sigma_n$; as before, once its size is fixed a key ring is formed by sampling a key 
pool with size $P_n$ randomly and without replacement. The graph $\mathbb{G}(n;\mathcal{D}_n,P_n)$ 
is then defined on the vertices $\{v_1,\ldots,v_n\}$ and contains an edge between any pair of nodes $v_x$ and $v_y$
as long as $\Sigma_x \cap \Sigma_y \neq \emptyset$. The model $\mathbb{G}(n;\boldsymbol{\mu}, \boldsymbol{K}_n,P_n)$ considered here constitutes a special case of $\mathbb{G}(n;\mathcal{D}_n,P_n)$ under the assumption that 
the support of $\mathcal{D}_n$ has a fixed size of $r$.

With these definitions in mind we now state the results by Zhao et al. \cite{ZhaoCDC} and by Goderhardt et al. \cite{GodehardtJaworskiRybarczyk}, respectively.

\begin{theorem}\cite[Theorem 1]{ZhaoCDC} \label{thm:zhao} {\sl Consider a general random intersection graph
$\mathbb{G}(n,\mathcal {D}_n,P_n)$. Let $|\Sigma|_n$ be a random variable following
the probability distribution $\mathcal {D}_n$. With a sequence $\alpha_n$
for all $n $ defined through\vspace{-1pt}
\begin{align}
 \frac{\mathbb{E}[|\Sigma|_n]^2}{P_n} & =
 \frac{\log  n + {(k-1)} \log \log n + {\alpha_n}}{n}, 
 \label{thm:grig:pe}
\end{align}
 if $\mathbb{E}[|\Sigma|_n] = \Omega\big(\sqrt{\log n}\hspace{2pt}\big)$, $\var{|\Sigma|_n} =
o\mathlarger{\mathlarger{\big(}}\frac{\mathbb{E}[|\Sigma|_n]^2}{ n(\log
n)^2 }\mathlarger{\mathlarger{\big)}}$ and $|\alpha_n| = o(\log n)$,
then
\begin{align}
 \lim_{n \to \infty}\mathbb{P} \big[\mathbb{G}(n,\mathcal {D}_n,P_n)\textrm{
is $k$-connected}\hspace{2pt}\big]  & =
\begin{cases} 0, &\textrm{ if $~\lim_{n \to \infty}{\alpha_n}
=-\infty$}, \\  1, &\textrm{ if $~\lim_{n \to \infty}{\alpha_n}
=\infty$,}  \\ e^{- \frac{e^{-\alpha ^*}}{(k-1)!}},
 &\textrm{ if $~\lim_{n \to \infty}{\alpha_n}
=\alpha ^* \in (-\infty, \infty)$.} \end{cases} 
\nonumber
 \end{align}
 }
 \end{theorem}
 
\begin{theorem}\cite[Theorem 2]{GodehardtJaworskiRybarczyk} \label{thm:godehardt} {\sl Consider a general random intersection graph
$\mathbb{G}(n,\mathcal {D},P_n)$, where $\mathcal{D}(\ell) = 0$ for all $\ell > r$ and $\ell=0$. Namely, all key ring
sizes are bound to be on the interval $[1, r]$. Let $|\Sigma|$ be a random variable following
the probability distribution $\mathcal {D}$. Then if
\begin{align}
\frac{n}{P_n}(\bE{|\Sigma|} -\mathcal{D}(1) ) - \log P_n \to \infty
\label{eq:scaling_godehardt}
\end{align}
then 
\[
\lim_{n \to \infty} \bP{\mathbb{G}(n,\mathcal {D},P_n)~\textrm{is connected}} = 1.
\]
Also, if 
$\mathcal{D}(r)=1$ for some $r \geq 2$, and it holds that
\begin{align}
n = P_n\frac{\log P_n+o(\log \log P_n)}{r^2},
\label{eq:scaling_godehardt_2}
\end{align}
then 
\[
\lim_{n \to \infty} \bP{\mathbb{G}(n,\mathcal {D},P_n)~\textrm{is connected}} = 0.
\]
  }
 \end{theorem}

In comparing Theorems \ref{thm:node_isolation}, \ref{thm:zhao} and \ref{thm:godehardt}, 
it is worth noting that $k$-connectivity is a stronger property than connectivity, which in turn is stronger than
absence of isolated nodes. However, although Theorems \ref{thm:zhao} and \ref{thm:godehardt}
consider strong graph properties, we now argue why the established results are {\em not} likely to be applicable 
for real-world sensor networks. First, Theorem \ref{thm:godehardt} focuses on the case where all 
possible key rings have a finite size that do not scale with $n$. In addition, with $\bE{|\Sigma|}$ fixed, it is clear that the
scaling conditions (\ref{eq:scaling_godehardt}) and (\ref{eq:scaling_godehardt_2}) both require
\begin{equation}
P_n = O\left(\frac{n}{\log n}\right).
\label{eq:key_pool_small}
\end{equation}
Unfortunately, it is often needed that key pool size $P_n$ be much larger than the network size $n$ \cite{virgil,DiPietroTissec} as otherwise the network will be extremely vulnerable against node capture attacks. 
In fact, one can see that with (\ref{eq:key_pool_small}) in effect, an adversary can compromise a significant portion of the key pool (and, hence network communication) by capturing $o(n)$ nodes.

We now focus on Theorem \ref{thm:zhao}, where the major problem arises from the assumption 
\begin{equation}
\var{|\Sigma|_n} =
o\left(\frac{\mathbb{E}[|\Sigma|_n]^2}{ n(\log
n)^2 }\right).
\label{eq:var_condition_zhao}
\end{equation}
For the model to be deemed as {\em inhomogeneous} random key graph, the variance of the key ring size
should be non-zero. In fact, given that key ring sizes are integer-valued, the simplest possible case would be that
$\mathcal{D}(K+1)=\mu$ and $\mathcal{D}(K)=1-\mu$ for some $0<\mu<1$ and positive integer $K$. This would amount to assigning either $K+1$ or $K$ keys to each node with probabilities $\mu$ and $1-\mu$, respectively. 
In this case, we can easily see that $\var{|\Sigma|} = \mu(1-\mu) > 0$ as long as $0<\mu<1$. Therefore, for an inhomogeneous random key graph, the condition (\ref{eq:var_condition_zhao}) implies that
$
\frac{\mathbb{E}[|\Sigma|_n]^2}{ n(\log
n)^2 } = w(1)$, 
 or, equivalently that
 \begin{equation}
 \bE{|\Sigma|_n} = w \left(\sqrt{n} \log n\right).
\label{eq:key_ring_large}
 \end{equation}
Put differently, Theorem \ref{thm:zhao} enforces {\em mean} key ring size to be much larger than $\sqrt{n} \log n$. 
However, a typical wireless sensor network will consist of a very large number of sensors, each with very limited memory and computational capability \cite{virgil,DiPietroTissec}. As a result, key rings with size 
$w(\sqrt{n} \log n)$ are unlikely to be implementable in most practical network deployments. In fact, it was suggested by Di Pietro et al. \cite{DiPietroTissec} that key rings with size $O(\log n)$ are acceptable for sensor networks. 

In comparison, our results Theorem \ref{thm:node_isolation} and Theorem \ref{thm:connectivity} 
do not require either of the unrealistic conditions (\ref{eq:key_pool_small}) or (\ref{eq:key_ring_large}). 
To see this, note that the enforced scaling condition (\ref{eq:scaling_law_2}) or its equivalent (\ref{eq:scaling_law_cor}) implies (see also Lemma \ref{lem:scaling_useful_cons}))
\begin{equation} 
\frac{K_{1,n} K_{r,n}}{P_n} =  \Theta\left(\frac{\log n}{n} \right).
\label{eq:condition_ultimate}
\end{equation}
It is clear that this condition does not require (\ref{eq:key_pool_small}), and in fact already enforces
$P_n = \Omega (n/\log n)$. As mentioned earlier, in real-world implementations the key pool size
is expected to grow at least linearly with $n$ so that the additional condition
 (\ref{eq:condition_for_con_1}) of Theorem \ref{thm:connectivity} is automatically satisfied. 
The second additional condition (\ref{eq:condition_for_con_2}) of 
our connectivity result and (\ref{eq:condition_ultimate}) can also be satisfied simultaneously
 without requiring the prohibitively large key ring sizes given at (\ref{eq:key_ring_large}).
To provide concrete examples, we can use $P_n = \Theta(n \log n)$, $K_{1,n}=\Theta(\log n)$ and $K_{r,n} = \Theta(\log n)$, or $P_n = \Theta(n \log n)$, $K_{1,n}=\Theta(\sqrt{\log n})$ and $K_{r,n} = \Theta((\log n)^{3/2})$. With proper choice of constants in these scalings, we will ensure that
i)  the resulting WSN is connected almost surely;
ii) the key pool size is much larger than the network
so that the resulting WSN has good level of resiliency against node capture attacks; and iii) the maximum key ring size used in the network is on the order of the ranges $\log n$ or $(\log n)^{3/2}$ that are usually regarded as 
feasible \cite{virgil,DiPietroTissec}; these choices also lead to
a much smaller mean key ring size than that prescribed in (\ref{eq:key_ring_large}).

In conclusion,
we showed that our results enable parameter choices that are widely regarded as practical in real-world sensor networks, while previous results given in  \cite{ZhaoCDC} and \cite{GodehardtJaworskiRybarczyk} do not.

\section{Preliminaries} 
\label{sec:prelim}

In this section, we establish several preliminary results that will be used in the proof of 
Theorem \ref{thm:node_isolation}. The first result states that mean edge probabilities are ordered in the same
way with the key ring sizes.

\begin{proposition}
{\sl
For any scaling  $K_1,\ldots,K_r, P: \mathbb{N}_0 \to \mathbb{N}_0^{r+1}$, we have
\begin{equation}
\lambda_1(n) \leq \lambda_2(n) \leq \cdots \leq \lambda_r(n)
\label{eq:ordering_lamda}
\end{equation}
for each $n=2,3,\ldots$.
\label{prop:ordering_lamda}
}
\end{proposition}

\myproof In view of (\ref{eq:mean_edge_prob_n}), the desired result (\ref{eq:ordering_lamda}) will follow immediately if we show that $p_{ij}(n)$ is increasing in both $i$ and $j$. Fix $n=2,3,\ldots$ and recall that $K_i$ increases as $i$
increases. For any $i,j$ such that $K_i +K_j > P$ we see from (\ref{eq:edge_prob}) that $p_{ij}(n)=1$; otherwise if $K_i+K_j \leq P$ we have $p_{ij}(n)<1$. Thus, given that $K_i + K_j$ increases with both $i$ and $j$, it will be enough to show that $p_{ij}(n)$ increases with both $i$ and $j$ on the range where $K_i + K_j \leq P$. But, on that range, 
we have
\begin{align}
{{{P-K_i}\choose{K_j}}\over{{P}\choose{K_j}}} 
& = \frac{(P-K_i)!}{P!} \frac{(P-K_j)!}{(P-K_i-K_j)!}
 = \prod_{\ell=0}^{K_i-1} \left(1-\frac{K_j}{P-\ell}\right).
\label{eq:inter_for_prelim_1}
\end{align}
It is now immediate that ${{{P-K_i}\choose{K_j}}\over{{P}\choose{K_j}}} 
$ decreases with both $K_i$ and $K_j$, and hence with $i$ and $j$. Hence, $p_{ij}(n)$ is seen to be increasing
with $i$ and $j$, and this establishes Proposition \ref{prop:ordering_lamda}. \myendpf

A useful consequence of Proposition \ref{prop:ordering_lamda} is given next.

\begin{lemma}
{\sl Consider any scaling  $K_1,\ldots,K_r, P: \mathbb{N}_0 \to \mathbb{N}_0^{r+1}$.
For any $i,j=1,\ldots, r$,
it holds that
\begin{equation}\nonumber
\lim_{n \rightarrow \infty} p_{ij}(n) = 0 ~\quad \textrm{if and only if} ~\quad
\lim_{n \rightarrow \infty} \frac{K_{i,n}K_{j,n}}{P_n} = 0,
\label{eq:Condition2}
\end{equation}
and under either condition we have the asymptotic equivalence
\begin{equation}\nonumber
p_{ij}(n) \sim \frac{K_{i,n}K_{j,n}}{P_n}.
\label{eq:AsymptoticsEquivalence}
\end{equation}
} \label{lem:AsymptoticEquivalence}
\end{lemma}
 {{\vspace{1mm}    \bf Proof. \ }}          
Lemma \ref{lem:AsymptoticEquivalence} can easily be established
by following the same arguments used in  \cite[Lemma 7.3]{yagan} or
\cite[Lemma 7.4.4]{YaganThesis}, namely by applying crude bounds (upper and lower)
to the expression (\ref{eq:inter_for_prelim_1}). The details are ommitted here for brevity.
\myendpf

Next, we give a result that collects several useful consequences of the scaling 
condition (\ref{eq:scaling_law}) under (\ref{eq:condition_for_con_2}).

\begin{lemma}
{\sl 
Consider any scaling  $K_1,\ldots,K_r, P: \mathbb{N}_0 \to \mathbb{N}_0^{r+1}$
such that (\ref{eq:scaling_law}) holds for some $c>0$. We have
\begin{align}
\frac{K_{1,n}K_{r,n}}{P_n} = \Theta\left(\frac{\log n}{n}\right).
\label{eq:bound_K1Kr_P}
\end{align}
If in addition (\ref{eq:condition_for_con_2}) holds, we have
\begin{align}
p_{rr}(n) \sim \frac{K_{r,n}^2}{P_n} = o\left(\frac{(\log n)^2}{n}\right).
\label{eq:K_r_square}
\end{align}
\label{lem:scaling_useful_cons}
}
\end{lemma}
 {{\vspace{1mm}    \bf Proof. \ }}          
The scaling condition
(\ref{eq:scaling_law_2}) states that
\[
\lambda_1(n) = \sum_{j=1}^{r} \mu_j p_{1j}(n) = c_n \frac{\log n}{n}
\]
with $\lim_{n \to \infty}c_n = c>0$. From the proof of Proposition (\ref{prop:ordering_lamda})
we know that $p_{ij}(n)$ increases with $i$ and $j$ under any scaling. Thus, we readily obtain that
\begin{align}
c_n \frac{\log n}{n} \leq p_{1r}(n) \leq \frac{c_n}{\mu_r} \frac{\log n}{n}.
\label{eq:inter_for_prelim_2}
\end{align}
Since $\mu_r>0$ this gives $p_{1r}(n) =\Theta\left(\frac{\log n}{n}\right)=o(1)$, whence we 
get $\frac{K_{1,n}K_{r,n}}{P_n} \sim p_{1r}(n)$
from Lemma \ref{lem:AsymptoticEquivalence} and (\ref{eq:bound_K1Kr_P}) is readily established.
We also find it useful to state the more detailed bounds
\begin{align}
\frac{c}{2} \frac{\log n}{n} \leq \frac{K_{1,n}K_{r,n}}{P_n} \leq  \frac{2c}{\mu_r}\left(\frac{\log n}{n}\right),
\label{eq:bound_K1Kr_P_detailed}
\end{align}
easily seen to be valid for any $n=2,3,\ldots$ sufficiently large in view of (\ref{eq:inter_for_prelim_2}).

We now turn to establishing (\ref{eq:K_r_square}). Comparing (\ref{eq:bound_K1Kr_P}) with (\ref{eq:condition_for_con_2}), we get
\begin{align}
\frac{K_{r,n}}{K_{1,n}}   =  {\frac{K_{1,n}K_{r,n}}{P_n} \over \frac{K_{1,n}^2}{P_n}} = \Theta \left(\frac{\log n}{w(1)}\right) = o(\log n).
\label{eq:bound_on_Kr_K1}
\end{align}
Next, we multiply (\ref{eq:bound_K1Kr_P}) with (\ref{eq:bound_on_Kr_K1}) to get
\[
\frac{K_{r,n}^2}{P_n}  =o \left(\frac{(\log n)^2}{n}\right) =o(1).
\]
Invoking Lemma \ref{lem:AsymptoticEquivalence} with $i=j=r$, we obtain (\ref{eq:K_r_square}).
\myendpf

The following inequality will also be useful in our proof.
\begin{proposition}
{\sl For any set of positive integers $K_1,\ldots,K_r$, $P$, and any
scalar $a \geq 1$, we have
\begin{align}
\frac{{P-\lceil a K_i \rceil \choose K_j}}{{P \choose K_j}} &\leq \left({{{P-K_i}\choose{K_j}}\over{{P}\choose{K_j}}}\right)^a, \quad i,j=1,\ldots, r.
 \label{eq:preliminary}
\end{align}
}
\label{prop:prelim2}
\end{proposition}
\myproof Fix $i,j =1,2,\ldots, r$. Observe that  ${{{P-K_i}\choose{K_j}}/{{P}\choose{K_j}}} \geq 0$ so that (\ref{eq:preliminary}) holds
trivially if $K_j + \lceil a K_i \rceil > P $. Assume here onwards that $K_j + \lceil a K_i \rceil \leq P$. Recalling 
(\ref{eq:inter_for_prelim_1}), we find
\begin{align}
\frac{{P-  \lceil a K_i \rceil \choose K_j}}{{P \choose K_j}} &=
\prod_{\ell=0}^{K_j-1} \left ( 1 - \frac{ \lceil a K_i \rceil}{P-\ell}
\right )
\leq \prod_{\ell=0}^{K_j-1} \left ( 1 - \frac{ a K_i }{P-\ell}
\right), 
\label{eq:prelim_1}
\end{align}
and
\begin{equation}
\frac{{P-   K_i  \choose K_j}}{{P \choose K_j}} = \prod_{\ell=0}^{K_j-1} \left ( 1 - \frac{ K_i }{P-\ell}
\right ). \label{eq:prelim_2}
\end{equation}

In view of (\ref{eq:prelim_1}) and (\ref{eq:prelim_2}), the
desired inequality (\ref{eq:preliminary}) will follow if we show
that
\begin{equation}
 1 - \frac{ a K_i }{P-\ell} \leq \left ( 1 - \frac{ K_i }{P-\ell}
\right ) ^ a, \quad \ell=0,1,\ldots, K_j-1. 
\label{eq:prelim_to_show}
\end{equation}
For each $\ell =0,1,\ldots, K_j-1$, (\ref{eq:prelim_to_show}) follows as we note  that
\[
1 - \left ( 1 - \frac{K_i}{P-\ell} \right )^{a} = \int_{ 1 -
\frac{K_i}{P-\ell} }^1 a t^{a-1}dt \leq \frac{a K_i}{P-\ell}
\]
and  (\ref{eq:preliminary}) is now established.
\myendpf

In the course of proving Theorem \ref{thm:node_isolation} we often make
use of the decomposition
\begin{equation}
\log ( 1 - x ) = -x - \Psi (x), \quad 0 \leq x < 1
\label{eq:LogDecomposition}
\end{equation}
with
$\Psi(x) := \int_0^x \frac{t}{1-t} dt$, and repeatedly 
use the fact that
\begin{equation}
\lim_{x \downarrow 0} \frac{ \Psi(x) }{x^2} = \frac{1}{2}.
\label{eq:limit_psi}
\end{equation}

Finally, we find it useful to derive a bound on  $p_{ij}$. Starting with (\ref{eq:inter_for_prelim_1})
we write
\begin{align}
\frac{{P-  K_i \choose K_j}}{{P \choose K_j}} = \prod_{\ell=0}^{K_i-1} \left(1-\frac{K_j}{P-\ell}\right)
\leq \prod_{\ell=0}^{K_i-1} \left(1-\frac{K_j}{P}\right) =\left(1-\frac{K_j}{P}\right)^{K_i} \leq e^{-\frac{K_i K_j}{P}}
\label{eq:easy_bound_on_q}
\end{align}

\section{A proof of Theorem \ref{thm:node_isolation} -- Establishing the zero-one law for absence of isolated nodes}
\label{sec:proof_thm_1}
The proof of Theorem \ref{thm:node_isolation} passes through applying the method of
first and second moments \cite[p.  55]{JansonLuczakRucinski}
to the number of isolated nodes 
in $\mathbb{G}(n; \boldsymbol{\mu}, \boldsymbol{K},P)$. To simplify the notation, we let $\boldsymbol{\theta}=(\boldsymbol{K},P)$. Let $I_n(\boldsymbol{\mu}, \boldsymbol{\theta})$ denote the total number of isolated nodes
in $\mathbb{G}(n; \boldsymbol{\mu},  \boldsymbol{\theta})$; i.e.,
\begin{equation}
I_n(\boldsymbol{\mu}, \boldsymbol{\theta}) = \sum_{\ell=1}^n \1{  v_\ell~{\rm is~isolated~in~} \mathbb{G}(n; \boldsymbol{\mu},  \boldsymbol{\theta}) }.
\label{eq:defn_I_n}
\end{equation}

\subsection{Establishing the one-law}
Consider now a scaling $\boldsymbol{\theta}: \mathbb{N}_0 \to \mathbb{N}_0^{r+1}$ such that 
(\ref{eq:scaling_law}) holds with $c>1$. 
The random graph $\mathbb{G}(n; \boldsymbol{\mu},  \boldsymbol{\theta}_n) $
 has no isolated
nodes if and only if $I_n(\boldsymbol{\mu}, \boldsymbol{\theta}_n)  = 0$.
The method of first moment \cite[Eqn. (3.10), p.
55]{JansonLuczakRucinski} gives
\begin{equation}
1 - \bE{ I_n(\boldsymbol{\mu}, \boldsymbol{\theta}_n)  } \leq \bP{  I_n(\boldsymbol{\mu}, \boldsymbol{\theta}_n)  = 0 },
\label{eq:FirstMoment}
\end{equation}
whence the one-law 
$
\lim_{n \to \infty} \bP{  I_n(\boldsymbol{\mu}, \boldsymbol{\theta}_n) = 0} =1
$
will follow if we show that
\begin{align}
\lim_{n \to \infty} \bE{ I_n(\boldsymbol{\mu}, \boldsymbol{\theta}_n) } = 0.
\label{eq:to_show_one_law}
\end{align}

By exchangeability of the indicator functions appearing at (\ref{eq:defn_I_n}), we find
\begin{align}
\bE{ I_n(\boldsymbol{\mu}, \boldsymbol{\theta}_n) } = n \bP{  v_1~{\rm is~isolated~in~} \mathbb{G}(n; \boldsymbol{\mu},  \boldsymbol{\theta}_n) }.
\label{eq:to_show_one_law_simplified}
\end{align}
Conditioning on the class of $v_1$, we further get
\begin{align}
n \bP{ v_1~{\rm is~isolated~in~} \mathbb{G}(n; \boldsymbol{\mu},  \boldsymbol{\theta}_n) }
&= n \sum_{i=1}^{r} \mu_i \bP{ v_1~{\rm is~isolated} ~ |~ v_1~{\rm is~class}~ i } 
\nonumber \\
&= n \sum_{i=1}^{r} \mu_i \bP{ \cap_{j=2}^{n}[v_1 \not \sim v_j]  ~ |~ v_1~{\rm is~class}~ i } 
\nonumber \\
&= n \sum_{i=1}^{r} \mu_i \left(\bP{v_1 \not \sim v_2  ~ |~ v_1~{\rm is~class}~ i }\right)^{n-1}
\label{eq:intermed_one_law} 
\end{align}
where (\ref{eq:intermed_one_law}) follows from the fact that rvs $\{v_1 \not \sim v_j\}_{j=2}^{n}$ are conditionally independent
given the key ring $\Sigma_1$ of node $v_1$. Conditioning further on the class of $v_2$, we find
\begin{align}
 \bP{v_1 \not \sim v_2  ~ |~ v_1~{\rm is~class}~ i }  
 &= \sum_{j=1}^{r} \mu_j \bP{v_1 \not \sim v_2  ~ |~ v_1~{\rm is~class}~ i,  v_2~{\rm is~class}~ j}
 \nonumber \\
 &= \sum_{j=1}^{r} \mu_j (1-p_{ij}(n))
  \nonumber \\
& = 1- \lambda_i(n).
\label{eq:inter2_one_law}
\end{align}

Using (\ref{eq:inter2_one_law}) in (\ref{eq:intermed_one_law}), and recalling
(\ref{eq:ordering_lamda})
we  get
\begin{align}\nonumber
 n \bP{ v_1~{\rm is~isolated~in~} \mathbb{G}(n; \boldsymbol{\mu},  \boldsymbol{\theta}_n) }
&= n \sum_{i=1}^{r} \mu_i (1- \lambda_i(n))^{n-1}
 \leq n (1- \lambda_1(n))^{n-1} 
  \leq e^{\log n - c_n \log n \frac{n-1}{n} }
\end{align}
as we also use (\ref{eq:scaling_law_2}). Letting $n$ go to infinity in this last expression we immediately get
\[
\lim_{n \to \infty} n \bP{ v_1~{\rm is~isolated~in~} \mathbb{G}(n; \boldsymbol{\mu},  \boldsymbol{\theta}_n) } 
= 0
\]
since 
$
\lim_{n \to \infty} 1- c_n \frac{n-1}{n} = 1-c < 0
$
under the enforced assumptions. Invoking (\ref{eq:to_show_one_law_simplified}) we now get (\ref{eq:to_show_one_law}) and the one-law is established.
\myendpf

\subsection{Establishing the zero-law}

This section is devoted to establishing the zero-law in Theorem \ref{thm:node_isolation}, namely the fact that
inhomogeneous random key graph contains at least one isolated node when the scaling condition (\ref{eq:scaling_law}) is satisfied with $c<1$. We will establish this 
by applying the method of second moment \cite[Remark 3.1, p.
55]{JansonLuczakRucinski} to a variable that counts nodes that are class-1 and isolated. Clearly, 
if we show that whp there exists at least one class-1 node that is isolated, then the desired zero-law will follow.   

Let $Y_n(\boldsymbol{\mu}, \boldsymbol{\theta})$ denote the number of isolated nodes in $\mathbb{G}(n; \boldsymbol{\mu},  \boldsymbol{\theta}_n) $ that are class-1. Namely, with $\chi_{n,i}(\boldsymbol{\mu}, \boldsymbol{\theta})$ denoting the indicator function that node $v_i$ is isolated and belongs to class-1, we have
$
Y_n(\boldsymbol{\mu}, \boldsymbol{\theta}) = \sum_{\ell=1}^{n} \chi_{n,\ell}(\boldsymbol{\mu}, \boldsymbol{\theta})$.  
The second moment method
states the inequality
\begin{equation}
\bP{  Y_n(\boldsymbol{\mu}, \boldsymbol{\theta})  = 0 } \leq 1 - \frac{ \bE{ Y_n(\boldsymbol{\mu}, \boldsymbol{\theta}) }^2}{
\bE{ Y_n(\boldsymbol{\mu}, \boldsymbol{\theta})  ^2} }. \label{eq:SecondMoment}
\end{equation}
Also, by exchangeability and the binary nature of the rvs $\chi_{n,1}(\boldsymbol{\mu}, \boldsymbol{\theta}), \ldots,
\chi_{n,n}(\boldsymbol{\mu}, \boldsymbol{\theta})$, we have
$\bE{ Y_n(\boldsymbol{\mu}, \boldsymbol{\theta})}= n \bE{\chi_{n,1}(\boldsymbol{\mu}, \boldsymbol{\theta}) }$
and
\begin{align}\label{eq:SecondMomentExpression}
\bE{ Y_n(\boldsymbol{\mu}, \boldsymbol{\theta}) ^2 } 
& = n \bE{\chi_{n,1}(\boldsymbol{\mu}, \boldsymbol{\theta}) }
 +  n(n-1) \bE{ \chi_{n,1}(\boldsymbol{\mu}, \boldsymbol{\theta})  \chi_{n,2}(\boldsymbol{\mu}, \boldsymbol{\theta}) }.
\end{align}
It then follows that
\begin{align}
\frac{ \bE{ Y_n(\boldsymbol{\mu}, \boldsymbol{\theta}) ^2 }}{ \bE{ Y_n(\boldsymbol{\mu}, \boldsymbol{\theta})  }^2 } &= \frac{
1}{ n\bE{\chi_{n,1}(\boldsymbol{\mu}, \boldsymbol{\theta})} }
  + \frac{n-1}{n} \cdot \frac{\bE{\chi_{n,1}(\boldsymbol{\mu}, \boldsymbol{\theta})
\chi_{n,2}(\boldsymbol{\mu}, \boldsymbol{\theta})}}
     {\left (  \bE{ \chi_{n,1}(\boldsymbol{\mu}, \boldsymbol{\theta})} \right )^2 }.
\label{eq:SecondMomentRatio}
\end{align}

From (\ref{eq:SecondMoment}) and
(\ref{eq:SecondMomentRatio}) we see that 
\begin{equation}
\lim_{n \to \infty} \bP{Y_n(\boldsymbol{\mu}, \boldsymbol{\theta}_n) = 0} = 0
\label{eq:intermed_zero_law}
\end{equation}
 holds if
\begin{equation}
\lim_{n \to \infty} n \bE{ \chi_{n,1}(\boldsymbol{\mu}, \boldsymbol{\theta}_n)}= \infty
\label{eq:OneLaw+NodeIsolation+convergence2}
\end{equation}
and
\begin{equation}
\limsup_{n \to \infty} \left( \frac{\bE{ \chi_{n,1}(\boldsymbol{\mu}, \boldsymbol{\theta}_n)
\chi_{n,2}(\boldsymbol{\mu}, \boldsymbol{\theta}_n)}}
     {\left (  \bE{ \chi_{n,1}(\boldsymbol{\mu}, \boldsymbol{\theta}_n) } \right )^2 }
\right) \leq 1. \label{eq:ZeroLaw+NodeIsolation+convergence}
\end{equation}
However, since $I_n(\boldsymbol{\mu}, \boldsymbol{\theta}_n)\geq Y_n(\boldsymbol{\mu}, \boldsymbol{\theta}_n)$,
(\ref{eq:intermed_zero_law}) immediately implies the
desired the zero-law 
\[
\lim_{n \to \infty} \bP{I_n(\boldsymbol{\mu}, \boldsymbol{\theta}_n) = 0} = 0.
\]

The next two technical propositions establish the needed results 
(\ref{eq:OneLaw+NodeIsolation+convergence2}) and
(\ref{eq:ZeroLaw+NodeIsolation+convergence}) under the appropriate
conditions on the scaling $\boldsymbol{\theta}: \mathbb{N}_0 \rightarrow
\mathbb{N}_0^{r+1}$.

\begin{proposition}
{\sl Consider a scaling  $K_1,\ldots,K_r, P: \mathbb{N}_0 \to \mathbb{N}_0^{r+1}$ such that 
(\ref{eq:scaling_law_2})
holds with $\lim_{n \to \infty} c_n = c>0$. Then, we have
\begin{equation}
n \bE{ \chi_{n,1}(\boldsymbol{\mu}, \boldsymbol{\theta}_n)} = (1+o(1)) \mu_1 n^{1-c_n}
\label{eq:bound_on_first_moment}
\end{equation}
so that
\begin{equation}
\lim_{n \rightarrow \infty } n \bE{ \chi_{n,1}(\boldsymbol{\mu}, \boldsymbol{\theta}_n)}  = \infty
\quad \textrm{if $~c<1$.} 
\label{eq:NodeIsolation+FirstMoment}
\end{equation}
 } \label{prop:Technical1}
\end{proposition}

\begin{proposition}
{\sl Consider a scaling  $K_1,\ldots,K_r, P: \mathbb{N}_0 \to \mathbb{N}_0^{r+1}$ such that 
(\ref{eq:scaling_law_2})
holds with $\lim_{n \to \infty} c_n = c>0$. Then, we have
(\ref{eq:ZeroLaw+NodeIsolation+convergence}).} \label{prop:Technical2}
\end{proposition}

A Proof of Proposition \ref{prop:Technical1} is given in Section
\ref{sec:ProofPropositionTechnical1}, while
Proposition \ref{prop:Technical2} is established in Section
\ref{sec:ProofPropositionTechnical2}.

\section{Proofs of Propositions \ref{prop:Technical1} and \ref{prop:Technical2}}
\label{sec:prop_proofs}
\subsection{A proof of Proposition \ref{prop:Technical1}}
\label{sec:ProofPropositionTechnical1}
Fix $n=2,3,\ldots$, and pick $\boldsymbol{u}$ and $\boldsymbol{\theta}$.
We have
\begin{align}
 n \bE{ \chi_{n,1}(\boldsymbol{\mu}, \boldsymbol{\theta})} 
  = n \bP{\textrm{$v_1$ is isolated and class-1}}
 &= n \mu_1  \bP{\textrm{$v_1$ is isolated} ~|~ \textrm{$v_1$ is class-1}}
\nonumber \\ 
& = n \mu_1 \bP{\cap_{j=2}^{n} [v_1 \not \sim v_{j}] ~|~ \textrm{$v_1$ is class-1}}
\nonumber \\ 
& = n \mu_1 \left( \bP{v_1 \not \sim v_{2}~|~ \textrm{$v_1$ is class-1}}\right)^{n-1}
\end{align}
by virtue of the fact that the events $\{v_1 \not \sim v_{j}\}_{j=2}^{n}$ are independent conditionally
on $\Sigma_1$. Invoking (\ref{eq:inter2_one_law}), we then get
\begin{equation}
n \bE{ \chi_{n,1}(\boldsymbol{\mu}, \boldsymbol{\theta})}  =  n \mu_1 \left(1-\lambda_1 \right)^{n-1}.
\label{eq:inter_1_prop_1}
\end{equation}

Now, consider a scaling  $K_1,\ldots,K_r, P: \mathbb{N}_0 \to \mathbb{N}_0^{r+1}$ such that 
(\ref{eq:scaling_law_2})
holds with $\lim_{n \to \infty} c_n = c>0$. Using this scaling in (\ref{eq:inter_1_prop_1}) and recalling (\ref{eq:LogDecomposition}) we get
\begin{align}
 n \bE{ \chi_{n,1}(\boldsymbol{\mu}, \boldsymbol{\theta}_n)} 
 = n \mu_1 \left(1-c_n \frac{\log n}{n} \right)^{n-1}
& = n \mu_1 e^{-c_n \log n \frac{n-1}{n} - (n-1)\psi\left(c_n \frac{\log n}{n}\right) }
\nonumber \\
& = \mu_1  n^{1-c_n} e^{ c_n \frac{\log n}{n}} e^{- (n-1) c_n^2 \frac{(\log n)^2}{n^2}\left(\frac{\psi\left(c_n \frac{\log n}{n}\right)}{\left(c_n \frac{\log n}{n}\right)^2} \right)}.
\end{align}
The desired result (\ref{eq:bound_on_first_moment}) is now immediate as we recall (\ref{eq:limit_psi})
and note that 
\[
\lim_{n \to \infty} c_n \frac{\log n}{n} = 0  {\textrm{,}} \quad \lim_{n \to \infty}\left(\frac{\psi\left(c_n \frac{\log n}{n}\right)}{\left(c_n \frac{\log n}{n}\right)^2} \right) =\frac{1}{2},
\quad \textrm{and}
\quad
\lim_{n \to \infty}(n-1) c_n^2 \frac{(\log n)^2}{n^2} = 0 
\]
since $\lim_{n \to \infty} c_n = c>0$.
From (\ref{eq:bound_on_first_moment}), we readily get (\ref{eq:NodeIsolation+FirstMoment}) upon noting that $\mu_1>0$.
\myendpf

\subsection{A proof of Proposition \ref{prop:Technical2}}
\label{sec:ProofPropositionTechnical2}

We start by obtaining an expression for the probability that nodes $v_1$ and $v_2$ are class-1 and isolated in  
$\mathbb{G}(n; \boldsymbol{\mu},  \boldsymbol{\theta})$. We get
\begin{align}
& \bE{ \chi_{n,1}(\boldsymbol{\mu}, \boldsymbol{\theta})
\chi_{n,2}(\boldsymbol{\mu}, \boldsymbol{\theta})}
\nonumber \\
& = \mu_1^2 \bP{\textrm{$v_1$ and $v_2$ are isolated} ~|~ \textrm{$v_1$ and $v_2$ are class-1}}
\nonumber \\
& = \mu_1^2 \bP{\Sigma_1 \cap \Sigma_2 = \emptyset ~\big|~ |\Sigma_1|=|\Sigma_2|=K_1}
  \bP{\cap_{j=3}^n [\Sigma_j \cap (\Sigma_1 \cup \Sigma_2) = \emptyset] 
~\bigg |~ 
\begin{array}{l}
  \Sigma_1 \cap \Sigma_2 = \emptyset, \\
  |\Sigma_1|=|\Sigma_2|=K_1  
\end{array}
}
\nonumber \\ 
& = \mu_1^2 {{{P-K_1}\choose{K_1}}\over{{P}\choose{K_1}}}
  \bP{\Sigma_3 \cap (\Sigma_1 \cup \Sigma_2) = \emptyset
~\bigg |
\begin{array}{l}
  \Sigma_1 \cap \Sigma_2 = \emptyset, \\
  |\Sigma_1|=|\Sigma_2|=K_1  
\end{array}
}^{n-2}
\nonumber \\
& = \mu_1^2 {{{P-K_1}\choose{K_1}}\over{{P}\choose{K_1}}} \left(\sum_{j=1}^{r} \mu_j {{{P-2K_1}\choose{K_j}}\over{{P}\choose{K_j}}}\right)^{n-2}
\label{eq:nodes_1_2_isolated}
\end{align} 
upon conditioning on the class of $v_3$. Similarly, it is easy to see that
\begin{equation}
\bE{ \chi_{n,1}(\boldsymbol{\mu}, \boldsymbol{\theta})} = \mu_1 \left(\sum_{j=1}^{r} \mu_j {{{P-K_1}\choose{K_j}}\over{{P}\choose{K_j}}}\right)^{n-1}.
\label{eq:inter1_prop2}
\end{equation}

Combining (\ref{eq:nodes_1_2_isolated}) and (\ref{eq:inter1_prop2}), we find
\begin{align}
\frac{\bE{ \chi_{n,1}(\boldsymbol{\mu}, \boldsymbol{\theta})
\chi_{n,2}(\boldsymbol{\mu}, \boldsymbol{\theta})}}
     {\left (  \bE{ \chi_{n,1}(\boldsymbol{\mu}, \boldsymbol{\theta}) } \right )^2 }
& = {{{P-K_1}\choose{K_1}}\over{{P}\choose{K_1}}}\left(\frac{\sum_{j=1}^{r} \mu_j {{{P-2K_1}\choose{K_j}}\over{{P}\choose{K_j}}}}{\left(\sum_{j=1}^{r} \mu_j {{{P-K_1}\choose{K_j}}\over{{P}\choose{K_j}}}\right)^2}\right)^{n-2}
\hspace{-3mm} \left(\sum_{j=1}^{r} \mu_j {{{P-K_1}\choose{K_j}}\over{{P}\choose{K_j}}}\right)^{-2}. 
\label{eq:key_expression_zero_law} 
\end{align}

Consider a scaling $\boldsymbol{\theta}: \mathbb{N}_0 \to \mathbb{N}_0^{r+1}$ such that 
(\ref{eq:scaling_law}) holds with $c < 1$. Reporting this scaling into the last expression, we see that 
\begin{align}
\left(\sum_{j=1}^{r} \mu_j {{{P_n-K_{1,n}}\choose{K_{j,n}}}\over{{P_n}\choose{K_{j,n}}}}\right)^{-2} 
\hspace{-2mm} = \left(1-\lambda_1(n) \right)^{-2}
  = \left(1-c_n \frac{\log n}{n}\right)^{-2} 
  = 1+o(1).
 \label{eq:zero_law_extra_term1}
\end{align}
With $p_{ij}(n)$ increasing with $i$ and $j$ as shown in Proposition \ref{prop:ordering_lamda}, it is also clear that
\begin{align}
 1 \geq {{{P_n-K_{1,n}}\choose{K_{1,n}}}\over{{P_n}\choose{K_{1,n}}}} = 1-p_{11}(n) \geq 1-\lambda_1(n)
 = 1-c_n \frac{\log n}{n},
 \nonumber
\end{align}
leading to
\begin{equation}
{{{P_n-K_{1,n}}\choose{K_{1,n}}}\over{{P_n}\choose{K_{1,n}}}}  = 1-o(1).
\label{eq:zero_law_extra_term2}
\end{equation}
Finally, we note from Proposition \ref{prop:prelim2} that
\begin{equation}
{{{P_n-2K_{1,n}}\choose{K_{j,n}}}\over{{P_n}\choose{K_{j,n}}}} \leq \left({{{P_n-K_{1,n}}\choose{K_{j,n}}}\over{{P_n}\choose{K_{j,n}}}}\right)^2, \qquad j=1,\ldots, r.
\label{eq:zero_law_extra_term3}
\end{equation}

Let $Z_n(\boldsymbol{\mu},\boldsymbol{\theta}_n)$ denote a rv such that
\[
Z_n(\boldsymbol{\mu},\boldsymbol{\theta}_n) = {{{P_n-K_{1,n}}\choose{K_{j,n}}}\over{{P_n}\choose{K_{j,n}}}} ~~~ \textrm{with probability $\mu_j$}, ~~ j=1,\ldots, r.
\]
Applying (\ref{eq:zero_law_extra_term1}), (\ref{eq:zero_law_extra_term2}), and (\ref{eq:zero_law_extra_term3})
in (\ref{eq:key_expression_zero_law}) we see that the desired result (\ref{eq:ZeroLaw+NodeIsolation+convergence}) will follow upon showing
\begin{equation}
\limsup_{n \to \infty}  \left(\frac{\bE{Z_n(\boldsymbol{\mu},\boldsymbol{\theta}_n)^2}}{\bE{Z_n(\boldsymbol{\mu},\boldsymbol{\theta}_n)}^2}\right)^{n-2} \leq 1. 
\label{eq:to_show_zero_law}
\end{equation}
We note that 
\begin{align}
\left(\frac{\bE{Z_n(\boldsymbol{\mu},\boldsymbol{\theta}_n)^2}}{\bE{Z_n(\boldsymbol{\mu},\boldsymbol{\theta}_n)}^2}\right)^{n-2}
 = \left(1+ \frac{\var{Z_n(\boldsymbol{\mu},\boldsymbol{\theta}_n)}}{\bE{Z_n(\boldsymbol{\mu},\boldsymbol{\theta}_n)}^2}\right)^{n-2}
& \leq e^{\frac{\var{Z_n(\boldsymbol{\mu},\boldsymbol{\theta}_n)} }{\bE{Z_n(\boldsymbol{\mu},\boldsymbol{\theta}_n)}^2}(n-2)}
\end{align}
 and that 
 \[
 \bE{Z_n(\boldsymbol{\mu},\boldsymbol{\theta}_n)} = 1-\lambda_1(n) =1-o(1).
 \]
 Hence, we will obtain (\ref{eq:to_show_zero_law}) if we show that
 \begin{equation}
 \lim_{n \to \infty} n \cdot \var{Z_n(\boldsymbol{\mu},\boldsymbol{\theta}_n)} = 0.
\label{eq:to_show_zero_law_simplified}
 \end{equation}
 
 In order to bound the variance of $Z_n(\boldsymbol{\mu},\boldsymbol{\theta}_n)$, we use Popoviciu's inequality
 \cite[p. 9]{jensen1999laguerre}. Namely, for any bounded rv $X$ with maximum value of $M$ and minimum value of $m$, we have
 \[
 \var{X} \leq \frac{1}{4}{(M-m)^2}.
 \]
 It is clear from the discussion given in the proof of Proposition \ref{prop:ordering_lamda} that 
 \[
{{{P_n-K_{1,n}}\choose{K_{r,n}}}\over{{P_n}\choose{K_{r,n}}}}  \leq Z_n(\boldsymbol{\mu},\boldsymbol{\theta}_n) \leq {{{P_n-K_{1,n}}\choose{K_{1,n}}}\over{{P_n}\choose{K_{1,n}}}} 
 \]
 holds for any scaling. Applying Popoviciu's inequality, we then get
 \begin{align}
 \var{Z_n(\boldsymbol{\mu},\boldsymbol{\theta}_n)}  \leq \frac{1}{4} \left({{{P_n-K_{1,n}}\choose{K_{1,n}}}\over{{P_n}\choose{K_{1,n}}}} - {{{P_n-K_{1,n}}\choose{K_{r,n}}}\over{{P_n}\choose{K_{r,n}}}} \right)^2
 & \leq \frac{1}{4} \left(1 - {{{P_n-K_{1,n}}\choose{K_{r,n}}}\over{{P_n}\choose{K_{r,n}}}} \right)^2
 = \frac{1}{4} \left(p_{1r}(n)\right)^2.
 \label{eq:var_inequality}
 \end{align}
 
Reporting the upper bound in (\ref{eq:inter_for_prelim_2}) into (\ref{eq:var_inequality}) we now find
 \begin{align}
n  \cdot \var{Z_n(\boldsymbol{\mu},\boldsymbol{\theta}_n)} \leq \frac{n}{4} \left(\frac{c_n}{\mu_r} \frac{\log n}{n}\right)^2.
 \end{align}
 Letting $n$ go to infinity in this last expression, we immediately get (\ref{eq:to_show_zero_law_simplified}) as we
 note that $\mu_r>0$ and $\lim_{n \to \infty}c_n = c > 0$. This establishes
 (\ref{eq:to_show_zero_law}) and 
 the desired result (\ref{eq:ZeroLaw+NodeIsolation+convergence}) now follows. \myendpf 
 
 \section{A proof of Theorem \ref{thm:connectivity} -- Establishing the zero-one law for connectivity}
\label{sec:proof_thm_2}
The proof of Theorem \ref{thm:connectivity} is technically more involved than that of Theorem \ref{thm:node_isolation}. To that end, we outline the main arguments leading to the proof in this section and complete the remaining steps in several 
sections that follow.

Fix $n=2,3, \ldots $ and consider 
$\boldsymbol{u}$ and $\boldsymbol{\theta}=(\boldsymbol{K},P)$.
We define the event
\[
C_n(\boldsymbol{u},\boldsymbol{\theta}) := \left [ \mathbb{G}(n;\boldsymbol{\mu},\boldsymbol{\theta})
\mbox{~is~connected} \right ]
\]
and recall that 
\[
[I_n( \boldsymbol{u},\boldsymbol{\theta})=0] = \left [ \mathbb{G}(n;\boldsymbol{\mu},\boldsymbol{\theta})
\mbox{~contains~no~isolated~nodes} \right ].
\]
If the random graph $\mathbb{G}(n;\boldsymbol{\mu},\boldsymbol{\theta})
$ is connected,
then it does not contain any isolated node, whence $C_n(\boldsymbol{\mu},\boldsymbol{\theta})$
is a subset of $[I_n( \boldsymbol{u},\boldsymbol{\theta})=0]$, and the conclusions
\begin{equation}
\bP{ C_n(\boldsymbol{\mu},\boldsymbol{\theta}) } \leq \bP{I_n( \boldsymbol{u},\boldsymbol{\theta})=0 }
\label{eq:FromConnectivityToNodeIsolation1}
\end{equation}
and
\begin{equation}
\bP{ C_n(\boldsymbol{\mu},\boldsymbol{\theta})^c } = \bP{ C_n(\boldsymbol{\mu},\boldsymbol{\theta})^c \cap(I_n( \boldsymbol{u},\boldsymbol{\theta})=0) } +
\bP{ (I_n( \boldsymbol{u},\boldsymbol{\theta})=0)^c } \label{eq:FromConnectivityToNodeIsolation2}
\end{equation}
follow.
Taken together with Theorem \ref{thm:node_isolation}, the
relations (\ref{eq:FromConnectivityToNodeIsolation1}) and
(\ref{eq:FromConnectivityToNodeIsolation2}) pave the way to
proving Theorem \ref{thm:connectivity}. To see this, pick any
scaling $K_1,\ldots,K_r, P: \mathbb{N}_0 \to \mathbb{N}_0^{r+1}$ such that 
(\ref{eq:scaling_law}) holds
for some $c>0$. If $c < 1$, then $\lim_{n \rightarrow
\infty} \bP{I_n( \boldsymbol{u},\boldsymbol{\theta}_n)=0} = 0$ by the zero-law for the absence
of isolated nodes (see Theorem \ref{thm:node_isolation}), whence $\lim_{n \rightarrow \infty} \bP{
C_n(\boldsymbol{\mu},\boldsymbol{\theta}_n) } = 0$ with the help of
(\ref{eq:FromConnectivityToNodeIsolation1}). If $c> 1$, then
$\lim_{n \rightarrow \infty} \bP{ I_n( \boldsymbol{u},\boldsymbol{\theta}_n)=0 } = 1$ by the
one-law for the absence of isolated nodes, and the desired
conclusion $\lim_{n \rightarrow \infty } \bP{ C_n(\boldsymbol{\mu},\boldsymbol{\theta}_n) } = 1$
(or equivalently, $\lim_{n \rightarrow \infty } \bP{
C_n(\boldsymbol{\mu},\boldsymbol{\theta}_n)^c } = 0$) will follow via
(\ref{eq:FromConnectivityToNodeIsolation2}) if we show the
following:

\begin{proposition}
{\sl For any scaling  $K_1,\ldots,K_r, P: \mathbb{N}_0 \to \mathbb{N}_0^{r+1}$
 such that
(\ref{eq:scaling_law}) holds for some $c>1$,
 we have
\begin{equation}
\lim_{n \rightarrow \infty} \bP{ C_n(\boldsymbol{\mu},\boldsymbol{\theta}_n)^c \cap (I_n( \boldsymbol{u},\boldsymbol{\theta}_n)=0) } = 0 . \label{eq:OneLawAfterReduction}
\end{equation}
as long as the conditions
(\ref{eq:condition_for_con_1}) and (\ref{eq:condition_for_con_2}) are satisfied. }
\label{prop:OneLawAfterReduction}
\end{proposition}

In words, Proposition \ref{prop:OneLawAfterReduction} 
states that the probability of the inhomogeneous random key graph being
not connected despite having no isolated nodes diminishes asymptotically under the enforced 
assumptions. In fact, the asymptotic equivalence of graph connectivity
and absence of isolated nodes is a well-known phenomenon in many
classes of random graphs; e.g., ER graphs \cite{bollobas}, random key graphs \cite{yagan},
and intersection of random key graphs and ER graphs \cite{yagan_onoff}. 

The basic idea in establishing Proposition
\ref{prop:OneLawAfterReduction} is to find a sufficiently tight
upper bound on the probability in (\ref{eq:OneLawAfterReduction})
and then to show that this bound goes to zero as $n$ becomes
large. Our approach is in the same vein with the one used for proving the
one-law for connectivity in ER graphs \cite[p. 164]{bollobas}. This 
approach has already proved useful in establishing one-laws for connectivity in
the standard random key graph \cite{yagan} and its intersection with an ER graph \cite{yagan_onoff}. 
Throughout the proof of the one-law for connectivity in {\em inhomogeneous} random key graphs, several
intermediate results will be borrowed directly from \cite{yagan,yagan_onoff} to avoid duplication.

We begin by deriving the needed upper bound on the term (\ref{eq:OneLawAfterReduction}).
Fix $n=2,3, \ldots $ and consider 
$\boldsymbol{u}$ and $\boldsymbol{\theta}=(\boldsymbol{K},P)$. For reasons that
will later become apparent we will need bounds on the number of {\em distinct} keys held by 
a specific set $S$ of sensors; just to give a hint, this will help us efficiently bound the probability that 
the sensors in $S$ are isolated from the rest of the network.
To that end, we define the
event $E_n(\boldsymbol{\mu},\boldsymbol{\theta};\boldsymbol{X})$ via
\begin{equation}
E_n(\boldsymbol{\mu},\boldsymbol{\theta};\boldsymbol{X})= \bigcup_{S \subseteq \mathcal{N}: ~
|S| \geq 1} ~ \left[\left|\cup_{i \in S}
\Sigma_i\right|~\leq~{X}_{|S|}\right]
\label{eq:E_n_defn}
\end{equation}
where $\mathcal{N}=\{1,\ldots,n\}$ and
$\boldsymbol{X}=[{X}_{1}~~{X}_{2}~~
\cdots~~ {X}_{n}]$ is an $n$-dimensional integer-valued
array. 

Let
\[
L_n := \min \left ( \left
\lfloor \frac{P}{K_1} \right \rfloor, \left \lfloor \frac{n}{2}
\right \rfloor \right )
\]
and set
\begin{eqnarray}
X_{\ell}= \left \{
\begin{array}{ll}
\lfloor \beta \ell K_1 \rfloor & ~ \mbox{$\ell=1,2,\ldots, L_n$} \\
 & \\
\lfloor \gamma P \rfloor &~ \mbox{$\ell=L_n+1, \ldots, n$}
\end{array}
\right . \label{eq:X_S_theta}
\end{eqnarray}
for some $\beta, \gamma$ in $(0,\frac{1}{2})$ that will be
specified later. 
With this setting, $E_n(\boldsymbol{\mu},\boldsymbol{\theta};\boldsymbol{X})$
encodes the event that for at least one $\ell=1,2,\ldots, n$, the total number of {\em distinct} keys held by 
at least one set of $\ell$ sensors is less than $\beta \ell K_1 \1{\ell \leq L_n} + \gamma P \1{\ell > L_n}$. Below, we
will show that by a careful selection of $\beta$ and $\gamma$, we can have the complement of $E_n(\boldsymbol{\mu},\boldsymbol{\theta}_n;\boldsymbol{X}_n)$ take place whp under the enforced assumptions on the
scaling $\boldsymbol{\theta}:\mathbb{N}_0 \to \mathbb{N}_0^{r+1}$; i.e., for {\em all} $\ell=1,2,\ldots,n$,
the total number of keys held by {\em any} set of $\ell$ sensors will be at least $\beta \ell K_1 \1{\ell \leq L_n} + \gamma P \1{\ell > L_n}$.
The relevance of this 
is easily seen as we use a simple bounding argument to write
\begin{align}\nonumber
&\bP{ C_n(\boldsymbol{\mu},\boldsymbol{\theta})^c \cap (I_n( \boldsymbol{u},\boldsymbol{\theta})=0) }
  \\ \nonumber
 & \leq
\bP{E_n(\boldsymbol{\mu},\boldsymbol{\theta};\boldsymbol{X})} + \bP{ C_n(\boldsymbol{\mu},\boldsymbol{\theta})^c \cap
(I_n( \boldsymbol{u},\boldsymbol{\theta})=0) \cap E_n(\boldsymbol{\mu},\boldsymbol{\theta};\boldsymbol{X})^c }.
\end{align}
It is now clear that the proof of Proposition \ref{prop:OneLawAfterReduction} and hence
that of Theorem \ref{thm:connectivity}
will consist of establishing the following two results.

\begin{proposition}
{\sl Consider a scaling $K_1,\ldots,K_r, P: \mathbb{N}_0 \to \mathbb{N}_0^{r+1}$ 
 such that
(\ref{eq:scaling_law}) holds for some $c>1$,
  (\ref{eq:condition_for_con_1b}) is satisfied for some
$\sigma>0$, and (\ref{eq:condition_for_con_2}) holds. We have
\begin{equation}
\lim_{n \rightarrow \infty} \bP{E_n(\boldsymbol{\mu},\boldsymbol{\theta}_n;\boldsymbol{X}_n)} =
0. \label{eq:OneLawAfterReductionPart1}
\end{equation}
where $\boldsymbol{X}_n=[X_{1,n} ~ \cdots ~
X_{n,n}]$ is as specified in (\ref{eq:X_S_theta}) with
$\beta$ in $(0, \frac{1}{2})$ is selected small enough to
ensure
\begin{equation}
\max \left ( 2 \beta \sigma , \beta \left( \frac{e^2}{\sigma}
\right) ^{\frac{ \beta }{ 1 - 2 \beta } } \right ) < 1,
\label{eq:ConditionOnLambda}
\end{equation}
and $\gamma$ in $(0, \frac{1}{2})$ is selected so that
\begin{equation}
\max \left ( 2 \left ( \sqrt{\gamma} \left ( \frac{e}{ \gamma } \right
)^{\gamma} \right )^\sigma, \sqrt{\gamma} \left ( \frac{e}{ \gamma }
\right)^{\gamma} \right ) < 1 . \label{eq:ConditionOnMU+1}
\end{equation}
 }
\label{prop:OneLawAfterReductionPart1}
\end{proposition}
A proof of Proposition \ref{prop:OneLawAfterReductionPart1} can be
found in Section \ref{sec:OneLawAfterReductionPart1}. Note that
for any $\sigma
>0$, $\lim_{\beta \downarrow 0} \beta \left( \frac{e^2}{\sigma}
\right) ^{\frac{ \beta }{ 1 - 2 \beta } } =0 $ so that the
condition (\ref{eq:ConditionOnLambda}) can always be met by
suitably selecting $\beta > 0$ small enough. Also, we have
$\lim_{\gamma \downarrow 0} \left ( \frac{e}{ \gamma } \right)^{\gamma}
=1$, whence $\lim_{\gamma \downarrow 0} \sqrt{\gamma} \left ( \frac{e}{
\gamma } \right)^{\gamma} = 0$, and (\ref{eq:ConditionOnMU+1}) can be
made to hold for any $\sigma>0$ by taking $\gamma > 0$ sufficiently
small.

\begin{proposition}
{\sl Consider a scaling $K_1,\ldots,K_r, P: \mathbb{N}_0 \to \mathbb{N}_0^{r+1}$ 
 such that
(\ref{eq:scaling_law}) holds for some $c>1$,
  (\ref{eq:condition_for_con_1b}) is satisfied for some
$\sigma>0$, and (\ref{eq:condition_for_con_2}) holds. We have
\begin{equation}
\lim_{n \rightarrow \infty} \bP{ C_n(\boldsymbol{\mu},\boldsymbol{\theta}_n)^c \cap
(I_n( \boldsymbol{u},\boldsymbol{\theta}_n)=0) \cap E_n(\boldsymbol{\mu},\boldsymbol{\theta}_n;\boldsymbol{X}_n)^c }
 = 0 .
\label{eq:OneLawAfterReductionPart2}
\end{equation}
where $\boldsymbol{X}_n=[X_{1,n} ~ \cdots ~
X_{n,n}]$ is as specified in (\ref{eq:X_S_theta}) with
$\gamma$ in $(0, \frac{1}{2})$ is selected small enough to ensure
(\ref{eq:ConditionOnMU+1}) and $\beta \in (0,\frac{1}{2})$ is
selected such that (\ref{eq:ConditionOnLambda}) is satisfied.
\label{prop:OneLawAfterReductionPart2} }
\end{proposition}

The proof of Proposition \ref{prop:OneLawAfterReductionPart2} is outlined in Section \ref{sec:OneLawAfterReductionPart2} with several steps completed in the sections that follow.


\section{A proof of Proposition \ref{prop:OneLawAfterReductionPart1}}
\label{sec:OneLawAfterReductionPart1}

The proof of Proposition \ref{prop:OneLawAfterReductionPart1} will follow similar steps 
to \cite[ Proposition 7.2]{yagan_onoff}
and
rely heavily on the results obtained by the author in
\cite{YaganThesis}. Using a standard union bound we first write
\begin{align}\nonumber
\bP{ E_n(\boldsymbol{\mu},\boldsymbol{\theta};\boldsymbol{X}) }
 &\leq  \sum_{ S
\subseteq \mathcal{N}: 1 \leq |S| \leq n } \bP{ \left|\cup_{i \in
S} \Sigma_{i}\right| \leq {X}_{|S|} }  = \sum_{\ell=1}^{ n } \left ( \sum_{S \in
\mathcal{N}_{n,\ell} } \bP{ \left|\cup_{i \in S}
\Sigma_{i}\right| \leq {X}_{\ell}} \right )
\end{align}
where $\mathcal{N}_{n,\ell} $ denotes the collection of all subsets
of $\{ 1, \ldots , n \}$ with exactly $\ell$ elements. 
Let $U_{\ell}(\boldsymbol{\mu},\boldsymbol{\theta})$ be given by
\[
U_{\ell}(\boldsymbol{\mu},\boldsymbol{\theta}) = |\cup_{i=1}^{\ell}\Sigma_i|, \qquad \ell=1,2,\ldots, n.
\]
By using
exchangeability and the fact that $|\mathcal{N}_{n,\ell} | = {n
\choose \ell}$, we get
\begin{align}\nonumber
\bP{ E_n(\boldsymbol{\mu},\boldsymbol{\theta};\boldsymbol{X}) } 
& \leq \sum_{\ell=1}^{ n } {n
\choose \ell}  \bP{ U_{\ell}(\boldsymbol{\mu},\boldsymbol{\theta})\leq{X}_{\ell} }
\\ \label{eq:E_n_big_expression}
&= \sum_{\ell=1}^{ L_n} {n \choose \ell}  \bP{
U_{\ell}(\boldsymbol{\mu},\boldsymbol{\theta})\leq\lfloor \beta \ell K_1 \rfloor }
 + \sum_{\ell=L_n+1}^{ n } {n \choose \ell} \bP{
U_{\ell}(\boldsymbol{\mu},\boldsymbol{\theta})\leq \lfloor \gamma P \rfloor }.
\end{align}

Now, consider any scaling 
$K_1,\ldots,K_r, P: \mathbb{N}_0 \to \mathbb{N}_0^{r+1}$ and
recall the ordering (\ref{eq:scaling_def}) of the key ring sizes. For any $\ell=1,2,\ldots, n$ 
define $U_{\ell}(K_{1,n},P_n)$ as
\[
U_{\ell}(K_{1,n},P_n) =_{st} U_{\ell}(\{1,0,0\ldots, 0\},\boldsymbol{\theta}_n).
\]
In other words, $U_{\ell}(K_{1,n},P_n)$ stands for the rv 
that has the same distribution with $U_{\ell}(\boldsymbol{\mu},\boldsymbol{\theta}_n)$ when $\boldsymbol{\mu}$
is {\em degenerate} with $\mu_1=1$ and $\mu_j=0$ for all $j=2,\ldots,r$; i.e., when all key ring
sizes are equal to $K_1$. With this setting,
$U_{\ell}(K_{1,n},P_n)$ is equivalent to the rv defined similarly
for the standard random key graph in \cite{yagan,yagan_onoff,YaganThesis}, where it was
often denoted by $U_{\ell}(\theta_n)$ with $\theta_n=(K_n,P_n)$. Given (\ref{eq:scaling_def}), it is a simple matter to
check that 
\begin{align}
U_{\ell}(K_{1,n},P_n) \preceq U_{\ell}(\boldsymbol{\mu},\boldsymbol{\theta}_n), \qquad \boldsymbol{\mu}=(
\mu_1,\ldots,\mu_r)
\label{eq:stochastic_order}
\end{align}
with $\preceq$ denoting the usual stochastic ordering. This can be seen by an
easy coupling argument where 
all sensors first receive $K_{1,n}$ keys and then an additional $K_{\ell,n}-K_{1,n}$ keys
are assigned to each sensor independently with probability $\mu_{\ell}$. Since additionally distributed 
keys can only increase the variable $U_{\ell}$, we obtain (\ref{eq:stochastic_order}).

Reporting (\ref{eq:stochastic_order}) into (\ref{eq:E_n_big_expression}) we now
get
\begin{align}\nonumber
\bP{ E_n(\boldsymbol{\mu},\boldsymbol{\theta}_n;\boldsymbol{X}_n)} 
& \leq  
\sum_{\ell=1}^{ L_n} {n \choose \ell}  \bP{
U_{\ell}(K_{1,n},P_n)\leq\lfloor \beta \ell K_{1,n} \rfloor }
 + \sum_{\ell=L_n+1}^{ n } {n \choose \ell} \bP{
U_{\ell}(K_{1,n},P_n)\leq \lfloor \gamma P_n \rfloor }.
\end{align}
Assume now that the scaling under consideration satisfies
(\ref{eq:scaling_law}) for some $c>1$, (\ref{eq:condition_for_con_1b}) with $\sigma>0$,
and
(\ref{eq:condition_for_con_2}). It was
shown \cite[Proposition 7.4.14, p. 142]{YaganThesis} that for any
scaling $K_1,P: \mathbb{N}_0 \rightarrow \mathbb{N}_0 \times \mathbb{N}_0$ such that
(\ref{eq:condition_for_con_1b}) holds for some
$\sigma>0$, we have
\begin{eqnarray}
\lim_{n \to \infty} \sum_{\ell= L_n+1}^{n} {n \choose \ell} ~
\bP{ U_{\ell}(K_{1,n},P_n) \leq \lfloor \gamma P_n \rfloor} = 0
\label{eq:E_n_part2}
\end{eqnarray}
whenever $\gamma$ in $(0,\frac{1}{2})$ is selected so that
(\ref{eq:ConditionOnMU+1}) holds; see also \cite[Proposition 7.4.17, p. 152]{YaganThesis}.
Hence, the desired conclusion
(\ref{eq:OneLawAfterReductionPart1}) will follow if we show that
\begin{equation}
\lim_{n \to \infty} \sum_{\ell= 1}^{L_n} {n \choose \ell} ~
\bP{ U_{\ell}(K_{1,n},P_n) \leq \lfloor \beta \ell K_{1,n} \rfloor} = 0
\label{eq:E_n_part1}
\end{equation}
under the condition (\ref{eq:ConditionOnLambda}). However, it can be seen from 
 \cite[Proposition 7.4.13, p. 140]{YaganThesis}
 and \cite[Proposition 7.4.16, p. 146]{YaganThesis}
 that with $\beta $ in $(0, \frac{1}{2})$
small enough to ensure (\ref{eq:ConditionOnLambda}) we have
(\ref{eq:E_n_part1}) for any scaling $K_1,P: \mathbb{N}_0 \rightarrow \mathbb{N}_0 \times \mathbb{N}_0$
such that $K_{1,n}=w(1)$. With this last condition clearly ensured under 
(\ref{eq:condition_for_con_1b}) and (\ref{eq:condition_for_con_2}) 
we obtain (\ref{eq:E_n_part1}). The proof of Proposition
\ref{prop:OneLawAfterReductionPart1} is now completed. \myendpf


\section{A proof of Proposition \ref{prop:OneLawAfterReductionPart2}}
\label{sec:OneLawAfterReductionPart2}

We will now work towards establishing (\ref{eq:OneLawAfterReductionPart2}), namely
 showing that the probability of $\mathbb{G}(n;\boldsymbol{\mu}_n,\boldsymbol{\theta}_n)$ being not connected 
despite having no isolated nodes approaches zero as $n$ gets large under the event
$E_n(\boldsymbol{\mu},\boldsymbol{\theta}_n;\boldsymbol{X}_n)^c$.
Fix $n=2,3, \ldots $ and consider 
$\boldsymbol{u}$ and $\boldsymbol{\theta}=(\boldsymbol{K},P)$.
For any non-empty subset $S$ of nodes, i.e., $S \subseteq \mathcal{V}= \{v_1,
\ldots , v_n \}$, we define the graph $\mathbb{G}(n;\boldsymbol{\mu},\boldsymbol{\theta})(S)$ 
(with vertex set $S$) as the subgraph of 
$\mathbb{G}(n;\boldsymbol{\mu},\boldsymbol{\theta})$ 
restricted to the nodes in $S$. 
With each non-empty subset $S$ of nodes, we associate several
events of interest: Let $C_n (\boldsymbol{\mu},\boldsymbol{\theta}; S)$ denote the event that
the subgraph $\mathbb{G}(n;\boldsymbol{\mu},\boldsymbol{\theta})(S)$ is itself
connected. It is clear that $C_n (\boldsymbol{\mu},\boldsymbol{\theta}; S)$ is completely determined
by the rvs $\{ \Sigma_i, \ v_i \in S \}$. 
We say that $S$
is {\em isolated} in $\mathbb{G}(n;\boldsymbol{\mu},\boldsymbol{\theta})$ if there are
no edges between the nodes in
$S$ and the nodes in the complement $S^c = \{v_1, \ldots , v_n\} -
S$. 
Let $B_n (\boldsymbol{\mu},\boldsymbol{\theta}; S)$
denote the event that $S$ is {\em isolated} in $\mathbb{G}(n;\boldsymbol{\mu},\boldsymbol{\theta})$, i.e.,
\begin{align}
B_n (\boldsymbol{\mu},\boldsymbol{\theta};S)
&:=  \left [ \Sigma_{i} \cap \Sigma_{j}
 = \emptyset, \quad v_i \in S , \ v_j \in S^c \right ] .
\nonumber
\end{align}
Finally, we set
\[
A_n (\boldsymbol{\mu},\boldsymbol{\theta};S) := C_n (\boldsymbol{\mu},\boldsymbol{\theta};S) 
\cap B_n (\boldsymbol{\mu},\boldsymbol{\theta};S) .
\]

Our main argument towards establishing  (\ref{eq:OneLawAfterReductionPart2})
relies on the following key
observation: If $\mathbb{G}(n;\boldsymbol{\mu},\boldsymbol{\theta})$ is {\em not}
connected and yet has {\em no} isolated nodes, then there must
exist a subset $S$ of nodes with $|S| \geq 2$ such that 
$\mathbb{G}(n;\boldsymbol{\mu},\boldsymbol{\theta})(S)$ is connected while $S$ is isolated in
$\mathbb{G}(n;\boldsymbol{\mu},\boldsymbol{\theta})$. This is captured by the inclusion
\begin{equation}
C_n(\boldsymbol{\mu},\boldsymbol{\theta})^c \cap (I_n( \boldsymbol{u},\boldsymbol{\theta})=0)\ \subseteq  \bigcup_{S \subseteq
\mathcal{V}: ~ |S| \geq 2} ~ A_n (\boldsymbol{\mu},\boldsymbol{\theta};S) \label{eq:BasicIdea}
\end{equation}
It is also clear that this union
need only be taken over all subsets $S$ of $\{v_1, \ldots , v_n \}$
with $2 \leq |S| \leq \lfloor \frac{n}{2} \rfloor $.

We now apply a standard union bound argument to (\ref{eq:BasicIdea}) and get
\begin{eqnarray}\nonumber
\lefteqn{\bP{ C_n(\boldsymbol{\mu},\boldsymbol{\theta})^c \cap (I_n( \boldsymbol{u},\boldsymbol{\theta})=0) \cap
E_n(\boldsymbol{\mu},\boldsymbol{\theta};\boldsymbol{X})^c }}
\\ \nonumber
 &\leq & \sum_{ S \subseteq
\mathcal{V}:~ 2 \leq |S| \leq \lfloor \frac{n}{2} \rfloor } \bP{
A_n (\boldsymbol{\mu},\boldsymbol{\theta};S) \cap E_n(\boldsymbol{\mu},\boldsymbol{\theta};\boldsymbol{X})^c }
\nonumber \\
&=& \sum_{\ell=2}^{ \lfloor \frac{n} {2} \rfloor } \left ( \sum_{S
\in \mathcal{V}_{n,\ell} } \bP{ A_n (\boldsymbol{\mu},\boldsymbol{\theta};S) \cap
E_n(\boldsymbol{\mu},\boldsymbol{\theta};\boldsymbol{X})^c} \right )
\label{eq:BasicIdea+UnionBound}
\end{eqnarray}
where $\mathcal{V}_{n,\ell} $ denotes the collection of all subsets
of $\{v_1, \ldots , v_n \}$ with exactly $\ell$ elements.

For each $\ell=1, \ldots , n$, we simplify the notation by writing
$A_{n,\ell} (\boldsymbol{\mu},\boldsymbol{\theta}) := A_n (\boldsymbol{\mu},\boldsymbol{\theta} ; \{ v_1, \ldots , v_\ell \} )$,
$B_{n,\ell} (\boldsymbol{\mu},\boldsymbol{\theta}) := B_n (\boldsymbol{\mu},\boldsymbol{\theta} ; \{ v_1, \ldots , v_\ell \} )$ and
$C_{n,\ell}(\boldsymbol{\mu},\boldsymbol{\theta}) := C_n (\boldsymbol{\mu},\boldsymbol{\theta} ; \{ v_1, \ldots , v_\ell \} )$. With a
slight abuse of notation, we use $C_n(\boldsymbol{\mu},\boldsymbol{\theta})$ for $\ell=n$ as
defined before. Under the enforced assumptions, exchangeability
yields
\[
\bP{ A_n (\boldsymbol{\mu},\boldsymbol{\theta};S) } = \bP{ A_{n,\ell} (\boldsymbol{\mu},\boldsymbol{\theta}) }, \quad S \in
\mathcal{V}_{n,\ell}
\]
and the expression
\begin{align} 
\sum_{S \in \mathcal{V}_{n,\ell} } \bP{ A_n (\boldsymbol{\mu},\boldsymbol{\theta};S)
\cap E_n(\boldsymbol{\mu},\boldsymbol{\theta};\boldsymbol{X})^c }&
= {n \choose \ell} ~ \bP{ A_{n,\ell} (\boldsymbol{\mu},\boldsymbol{\theta}) \cap
E_n(\boldsymbol{\mu},\boldsymbol{\theta};\boldsymbol{X})^c } \label{eq:ForEach=r}
\end{align}
follows since $|\mathcal{V}_{n,\ell} | = {n \choose \ell}$. Substituting
into (\ref{eq:BasicIdea+UnionBound}) we obtain the key bound
\begin{align}
\bP{ C_n(\boldsymbol{\mu},\boldsymbol{\theta})^c \cap (I_n( \boldsymbol{u},\boldsymbol{\theta})=0) \cap
E_n(\boldsymbol{\mu},\boldsymbol{\theta};\boldsymbol{X})^c }
 &\leq \sum_{\ell=2}^{ \lfloor
\frac{n}{2} \rfloor } {n \choose \ell} ~ \bP{ A_{n,\ell}(\boldsymbol{\mu},\boldsymbol{\theta})  \cap
E_n(\boldsymbol{\mu},\boldsymbol{\theta};\boldsymbol{X})^c} .
\label{eq:BasicIdea+UnionBound2}
\end{align}


Next, we derive bounds for the probabilities appearing at (\ref{eq:BasicIdea+UnionBound2}).
Recall the definitions (\ref{eq:edge_prob}) and (\ref{eq:mean_edge_prob}).
\begin{proposition}
{\sl
Consider $\boldsymbol{\theta}=(\boldsymbol{K},P)$ and $\boldsymbol{\mu}=(\mu_1, \ldots, \mu_r)$
such that $K_1 \leq K_2 \leq \cdots \leq K_r$. We have
\begin{align}
\bP{ A_{n,\ell}(\boldsymbol{\mu},\boldsymbol{\theta} ) \cap
E_n(\boldsymbol{\mu},\boldsymbol{\theta};\boldsymbol{X})^c} \leq \min\left\{1, \ell^{\ell-2} (p_{rr})^{\ell-1}\right\}
\left(\min\left\{1-\lambda_1,\bE{e^{-\frac{(X_{\ell}+1)|\Sigma|}{P}}} \right\}\right)^{n-\ell}
\label{eq:new_combined_bounds}
\end{align}
where $|\Sigma|$ denotes a rv such that
\[
|\Sigma| = K_j \quad \textrm{with probability $\mu_j$}, \quad j=1,\ldots, r.
\]
}
\label{prop:new_combined_bounds}
\end{proposition}
 {{\vspace{1mm}    \bf Proof. \ }}          
We start by observing the equivalence
\[
B_{n,\ell}(\boldsymbol{\mu},\boldsymbol{\theta})= \left [ \left ( \cup_{i =1}^{\ell}
\Sigma_i \right ) \cap \Sigma_j = \emptyset, ~~\ j=\ell+1, \ldots,
n \right ].
\]
Hence,
under the enforced assumptions on the rvs $\Sigma_1, \ldots ,
\Sigma_n $, we readily obtain the expression
\begin{align}\nonumber
\bP{ B_{n,\ell}(\boldsymbol{\mu},\boldsymbol{\theta}) ~~\Big | ~~ 
  \Sigma_1, \ldots , \Sigma_{\ell}   }
&= \prod_{j=\ell+1}^n \bE{  {{P- |\cup_{i =1}^{\ell}
\Sigma_i|} \choose {|\Sigma_j|}} \over {{P \choose |\Sigma_j|}} } =
 \bE{ {{P- |\cup_{i =1}^{\ell}
\Sigma_i|} \choose {|\Sigma|}} \over {{P \choose |\Sigma|}} }^{n-\ell}
\end{align}
where $|\Sigma|$ denotes a rv that takes the value $K_j$ with probability $\mu_j$ for
each $j=1,\ldots,r$. Note that we always have 
$|\cup_{i =1}^{\ell}
\Sigma_i| \geq K_1$, while it holds that $|\cup_{i =1}^{\ell}
\Sigma_i| \geq X_{\ell}+1$ on the event $E_n(\boldsymbol{\mu},\boldsymbol{\theta};\boldsymbol{X})^c$. Combining,
we get
\begin{align}
\bP{ B_{n,\ell}(\boldsymbol{\mu},\boldsymbol{\theta}) \cap E_n(\boldsymbol{\mu},\boldsymbol{\theta};\boldsymbol{X})^c ~~\Big | ~~ 
  \Sigma_1, \ldots , \Sigma_{\ell}   }
&\leq  \left(\min\left\{\bE{  {{P- K_1} \choose {|\Sigma|}} \over {{P \choose |\Sigma|}} },\bE{  {{P- (X_{\ell}+1)} \choose {|\Sigma|}} \over {{P \choose |\Sigma|}} }\right\}\right)^{n-\ell}
\nonumber \\
& \leq \left(\min\left\{1-\lambda_1,\bE{e^{-\frac{(X_{\ell}+1)|\Sigma|}{P}}} \right\}\right)^{n-\ell},
\label{eq:U_r_Difficult_Bould}
\end{align}
where in the last step we also used (\ref{eq:easy_bound_on_q}).

Conditioning on the rvs $\Sigma_1,
\ldots , \Sigma_{\ell}$  which fully determine the event $C_{n,\ell}(\boldsymbol{\mu},\boldsymbol{\theta}))$, we conclude via
(\ref{eq:U_r_Difficult_Bould}) that
\begin{align}\nonumber
\bP{  A_{n,\ell}(\boldsymbol{\mu},\boldsymbol{\theta} ) \cap
E_n(\boldsymbol{\mu},\boldsymbol{\theta};\boldsymbol{X})^c}
  &= \bP{ C_{n,\ell}(\boldsymbol{\mu},\boldsymbol{\theta})) \cap B_{n,\ell}(\boldsymbol{\mu},\boldsymbol{\theta})) \cap
E_n(\boldsymbol{\mu},\boldsymbol{\theta};\boldsymbol{X})^c}
 \\ \nonumber
&= \bE{ \1{C_{n,\ell}(\boldsymbol{\mu},\boldsymbol{\theta}))} 
\bP{ B_{n,\ell}(\boldsymbol{\mu},\boldsymbol{\theta}) \cap E_n(\boldsymbol{\mu},\boldsymbol{\theta};\boldsymbol{X})^c ~~\Big | ~~ 
  \Sigma_1, \ldots , \Sigma_{\ell}   }
}
 \\ \nonumber
&\leq \bP{C_{n,\ell}(\boldsymbol{\mu},\boldsymbol{\theta}))}
\left(\min\left\{1-\lambda_1,\bE{e^{-\frac{(X_{\ell}+1)|\Sigma|}{P}}} \right\}\right)^{n-\ell}.
\end{align}
In view of this last bound, Proposition \ref{prop:new_combined_bounds} will be established if we show that
\begin{align}
\bP{C_{n,\ell}(\boldsymbol{\mu},\boldsymbol{\theta}))} \leq \ell^{\ell-2} \left(p_{rr}\right)^{\ell-1}.
\label{eq:upper_bound_C_n}
\end{align}

To  see why (\ref{eq:upper_bound_C_n}) holds, observe that the subgraph of 
$\mathbb{G}(n;\boldsymbol{\mu},\boldsymbol{\theta})$ on the vertices $v_1,\ldots,v_{\ell}$, 
hereafter denoted $\mathbb{G}_{\ell}(n;\boldsymbol{\mu},\boldsymbol{\theta})$,
is connected if and only if it contains a {\em spanning tree}. 
Let ${\cal
T}_{\ell}$ denote the collection of all spanning trees on the vertex
set $\{v_1, \ldots , v_\ell \}$.
Then, we have
\[
C_{n,\ell}(\boldsymbol{\mu},\boldsymbol{\theta}) = \cup_{T \in \cal
{T}_{\ell}}   [T \subseteq \mathbb{G}_{\ell}(n;\boldsymbol{\mu},\boldsymbol{\theta})]
\]
By Cayley's formula \cite{martin2013counting} there are $\ell^{\ell-2}$ trees on $\ell$
vertices, i.e., $| {\cal T}_{\ell}| = \ell^{\ell-2}$. In addition, for any $T$ in $\mathcal{T}_{\ell}$,
it is clear that 
\[
\bP{T \subseteq \mathbb{G}_{\ell}(n;\boldsymbol{\mu},\boldsymbol{\theta})} \leq 
\bP{T \subseteq \mathbb{G}_{\ell}(n;\{0,\ldots,0,1\},\boldsymbol{\theta})}
\]
since $K_r \geq K_j$ for any $j=1,2,\ldots,r-1$; i.e., probability that the tree $T$ is contained
in this subgraph is maximized when all the nodes in the subgraph belong to class-$r$ that are assigned
the largest number of keys. With $\mu_r=1$, $\mathbb{G}_{\ell}(n;\boldsymbol{\mu},\boldsymbol{\theta})$
becomes equivalent to the standard random key graph $\mathbb{G}_{\ell}(n;K_r,P)$ for 
which it is known \cite[Lemma 9.1]{yagan} that 
\[
\bP{T \subseteq \mathbb{G}_{\ell}(n;\boldsymbol{\mu},\boldsymbol{\theta})} = \left(p_{rr}\right)^{\ell-1},\qquad
T \in \mathcal{T}_{\ell}, ~~ \ell=2,3,\ldots
\]
This follow from the fact that a tree on $\ell$ nodes consists of $\ell-1$ edges and that edge events in the random
key graph are pairwise independent. Collecting, we obtain via a union bound that
\[
\bP{C_{n,\ell}(\boldsymbol{\mu},\boldsymbol{\theta})} \leq \sum_{T \in \cal
{T}_{\ell}}   \bP{T \subseteq \mathbb{G}_{\ell}(n;\boldsymbol{\mu},\boldsymbol{\theta})}
\leq \sum_{T \in \cal
{T}_{\ell}}  \bP{T \subseteq \mathbb{G}_{\ell}(n;\{0,\ldots,0,1\},\boldsymbol{\theta})} \leq \ell^{\ell-2}\left(p_{rr}\right)^{\ell-1}.
\]
This establishes (\ref{eq:upper_bound_C_n}) and the proof of
Proposition (\ref{prop:new_combined_bounds}) is now completed.\myendpf

Now, consider a scaling  $\boldsymbol{\theta}: \mathbb{N}_0 \to \mathbb{N}_0^{r+1}$ as in the statement of
Proposition \ref{prop:OneLawAfterReductionPart2}. Using this scaling in 
(\ref{eq:BasicIdea+UnionBound2}) together with (\ref{eq:new_combined_bounds}) we see that the proof of Proposition
\ref{prop:OneLawAfterReductionPart2}
will be completed once we show
\begin{equation}
\lim_{n \rightarrow \infty} \sum_{\ell=2}^{ \lfloor \frac{n}{2}
\rfloor } {n \choose \ell} \min\left\{1, \ell^{\ell-2} (p_{rr}(n))^{\ell-1}\right\}
\left(\min\left\{1-\lambda_1(n),\bE{e^{-\frac{(X_{\ell,n}+1)|\Sigma|_n}{P_n}}} \right\}\right)^{n-\ell}= 0. \label{eq:OneLawToShow}
\end{equation}
Combined with Proposition \ref{prop:OneLawAfterReductionPart1},
this will lead to Proposition \ref{prop:OneLawAfterReduction} and hence to Theorem \ref{thm:connectivity}.
To that end, we devote the rest of the paper to establishing (\ref{eq:OneLawToShow}). Throughout, we make repeated use of the standard
bounds
\begin{equation}
{n \choose \ell} \leq \left ( \frac{e n}{\ell} \right )^{\ell}, \quad
\begin{array}{c}
\ell=1, \ldots , n \\
n=1,2, \ldots \\
\end{array}
\label{eq:CombinatorialBound1}
\end{equation}
and
\begin{equation}
\sum_{\ell=2}^{\lfloor n/2 \rfloor}{n \choose \ell} \leq 2^n,
\label{eq:CombinatorialBound2}
\end{equation}
where the latter follows from the Binomial formula.


\section{Establishing (\ref{eq:OneLawToShow})}
\label{sec:OneLawAfterReductionPart3}

We will establish  (\ref{eq:OneLawToShow}) in several steps with each step focusing on a specific range 
of the summation over $\ell$. 
Throughout this section, we consider a scaling  $K_1,\ldots,K_r, P: \mathbb{N}_0 \to \mathbb{N}_0^{r+1}$ such that 
(\ref{eq:scaling_law}) holds for some $c>1$, (\ref{eq:condition_for_con_1}) holds for some $\sigma>0$ and
(\ref{eq:condition_for_con_2}) is satisfied. The desired result (\ref{eq:OneLawToShow}) 
follows from (\ref{eq:desired_limit_range_1}), (\ref{eq:desired_limit_range_2}), (\ref{eq:desired_limit_range_3}),
(\ref{eq:desired_limit_range_4}), and
(\ref{eq:desired_limit_range_5}) that are
established in Sections \ref{subsec:range1}-\ref{subsec:range5}, respectively.

\subsection{The case where $2 \leq \ell \leq R$}
\label{subsec:range1}
The first range considers fixed values of $\ell$.
For the moment, fix an integer $R$ that will be specified later in Section \ref{subsec:range2}; see (\ref{eq:choosing_R}).
For each $\ell=2,\ldots, R$ we use (\ref{eq:CombinatorialBound1}), (\ref{eq:K_r_square}), and (\ref{eq:scaling_law})
to get
\begin{align}
& {n \choose \ell} \min\left\{1, \ell^{\ell-2} (p_{rr}(n))^{\ell-1}\right\}
\left(\min\left\{1-\lambda_1(n),\bE{e^{-\frac{(X_{\ell,n}+1)|\Sigma|_n}{P_n}}} \right\}\right)^{n-\ell}
\nonumber \\
& ~ \leq \left ( \frac{e n}{\ell} \right )^{\ell} \ell^{\ell-2} (p_{rr}(n))^{\ell-1} \left(1-\lambda_1(n)\right)^{n-\ell}
\nonumber \\
& ~\leq  (e n)^{\ell}(p_{rr}(n))^{\ell-1}  e^{-c_n \log n \frac{n-\ell}{n}}
\nonumber \\
& ~= o(1) n^{\ell} \left(\frac{(\log n)^2}{n}\right)^{\ell-1} n^{-c_n  \frac{n-\ell}{n}}
\nonumber \\
& ~= o(1) n^{1-c_n  \frac{n-\ell}{n}} (\log n)^{2\ell-2}
\nonumber \\
& ~= o(1)
\end{align}
upon noting that
$
\lim_{n \to \infty} 1-c_n  \frac{n-\ell}{n} = 1-c < 0
$ on the given range of $\ell$. Thus, for any $R$ we have
\begin{align}
\lim_{n \rightarrow \infty} \sum_{\ell=2}^{ R } {n \choose \ell} \min\left\{1, \ell^{\ell-2} (p_{rr}(n))^{\ell-1}\right\}
\left(\min\left\{1-\lambda_1(n),\bE{e^{-\frac{(X_{\ell,n}+1)|\Sigma|_n}{P_n}}} \right\}\right)^{n-\ell}= 0.
\label{eq:desired_limit_range_1}
\end{align}

\subsection{The case where $R+1 \leq \ell \leq \left \lfloor \frac{\mu_r n}{2 \beta c \log n} \right \rfloor$}
\label{subsec:range2}
Next, we handle the range where $R+1 \leq \ell \leq  \left \lfloor \frac{\mu_r n}{2 \beta c \log n}\right \rfloor$. Noting that 
$K_{1,n} \leq K_{r,n}$, we realize from
(\ref{eq:bound_K1Kr_P}) that 
$\frac{K_{1,n}^2}{P_n} = O\left(\frac{\log n}{n}\right)$, or equivalently
that $\frac{P_n}{K_{1,n}^2} = \Omega \left(\frac{n}{\log n}\right)$. Also, (\ref{eq:condition_for_con_1}) and
(\ref{eq:condition_for_con_2}) imply together that $K_{1,n}=w(1)$, leading to
$\frac{P_n}{K_{1,n}} = w \left(\frac{n}{\log n}\right)$.
Hence, on the range under consideration here, we have $\ell \leq L_n= \min  ( 
\lfloor \frac{P_n}{K_{1,n}}  \rfloor,  \lfloor \frac{n}{2}
 \rfloor  )$
 so that $X_{\ell,n} = \lfloor \beta \ell K_{1,n} \rfloor$. With this in mind, we get
\begin{align}
\bE{e^{-\frac{(X_{\ell,n}+1)|\Sigma|_n}{P_n}}} \leq \bE{e^{-\frac{\beta \ell K_{1,n}|\Sigma|_n}{P_n}}} =
\sum_{j=1}^{r} \mu_j e^{-\frac{\beta \ell K_{1,n}K_{j,n}}{P_n}} \leq 1-\mu_r  
+ \mu_r e^{-\frac{\beta \ell K_{1,n}K_{r,n}}{P_n}}. 
\label{eq:bounding_expectation}
\end{align}
In view of (\ref{eq:bound_K1Kr_P_detailed}),
\[
\beta \ell  \frac{K_{1,n}K_{r,n}}{P_n} \leq \beta \ell \frac{2 c}{\mu_r} \frac{\log n}{n} \leq 1.
\]
holds for all $n$ sufficiently large and $\ell \leq \frac{\mu_r n}{2 \beta c \log n}$. Hence, on the same 
range we have 
\[
1-e^{-\frac{\beta \ell K_{1,n}K_{r,n}}{P_n}} \geq \frac{\beta \ell K_{1,n}K_{r,n}}{2P_n}. 
\]
Reporting these into (\ref{eq:bounding_expectation}) we get
\begin{align}
\bE{e^{-\frac{(X_{\ell,n}+1)|\Sigma|_n}{P_n}}} \leq 1-\mu_r\left(1-e^{-\frac{\beta \ell K_{1,n}K_{r,n}}{P_n}}  \right)
\leq 1-\mu_r  \frac{\beta \ell K_{1,n}K_{r,n}}{2P_n} \leq e^{-\frac{\mu_r  \beta c}{4} \ell  \frac{\log n}{n} } 
\label{eq:range_2_inter_1}
\end{align}
for all $n$ sufficiently large, where the last inequality follows from the lower bound in (\ref{eq:bound_K1Kr_P_detailed}). 
%

Consider now the range of $n$ sufficiently large that (\ref{eq:range_2_inter_1}) is valid. Using (\ref{eq:range_2_inter_1}), 
(\ref{eq:CombinatorialBound1}), and (\ref{eq:K_r_square}), we get
\begin{align}
& \sum_{\ell=R+1}^{ \left \lfloor\frac{\mu_r n}{2 \beta c \log n} \right \rfloor}{n \choose \ell} \min\left\{1, \ell^{\ell-2} (p_{rr}(n))^{\ell-1}\right\}
\left(\min\left\{1-\lambda_1(n),\bE{e^{-\frac{(X_{\ell,n}+1)|\Sigma|_n}{P_n}}} \right\}\right)^{n-\ell}
\nonumber \\
& ~ \leq \sum_{\ell=R+1}^{ \left \lfloor\frac{\mu_r n}{2 \beta c \log n} \right \rfloor} \left ( \frac{e n}{\ell} \right )^{\ell} \ell^{\ell-2}  \left(\frac{(\log n)^2}{n}\right)^{\ell-1}
e^{-\frac{\mu_r  \beta c}{4} \ell  \frac{\log n}{n} (n-\ell)}
\nonumber \\
& ~\leq  n \sum_{\ell=R+1}^{ \left \lfloor\frac{\mu_r n}{2 \beta c \log n} \right \rfloor} \left(e (\log n)^2 \right)^{\ell}  e^{-\frac{\mu_r  \beta c}{4} \ell  \frac{\log n}{n} (n-\lfloor \frac{n}{2}\rfloor)}
\nonumber \\
& ~\leq   n \sum_{\ell=R+1}^{ \infty}  \left(e (\log n)^2 e^{-\frac{\mu_r  \beta c}{8} \log n} \right)^{\ell}. 
\label{eq:range_2_inter_2}  
\end{align}
Since $\mu_r$, $\beta$, $c$ are all positive scalars, we have
\[
\lim_{n \to \infty}e (\log n)^2 e^{-\frac{\mu_r  \beta c}{8} \log n} = 0
\]
so that the infinite series appearing at (\ref{eq:range_2_inter_2}) is summable. In fact, we have
\begin{align}
 n \sum_{\ell=R+1}^{ \infty}  \left(e (\log n)^2 e^{-\frac{\mu_r  \beta c}{8} \log n} \right)^{\ell} &= (1+o(1)) n \left(e (\log n)^2 e^{-\frac{\mu_r  \beta c}{8} \log n} \right)^{R+1} 
 \nonumber \\
 & = O(1) n^{1-(R+1)\frac{\mu_r  \beta c}{8}}  (\log n)^{2R+2}.
\end{align}
It is now clear that we get
\begin{align}
\lim_{n \to \infty}  \sum_{\ell=R+1}^{ \left \lfloor\frac{\mu_r n}{2 \beta c \log n} \right \rfloor}{n \choose \ell} \min\left\{1, \ell^{\ell-2} (p_{rr}(n))^{\ell-1}\right\}
\left(\min\left\{1-\lambda_1(n),\bE{e^{-\frac{(X_{\ell,n}+1)|\Sigma|_n}{P_n}}} \right\}\right)^{n-\ell}
 = 0
 \label{eq:desired_limit_range_2}
\end{align}
as long as $R$ satisfies
\begin{align}
R \geq \left\lceil \frac{8}{\mu_r  \beta c}\right\rceil.
\label{eq:choosing_R}
\end{align}
Since $\mu_r, \beta, c > 0$, such a selection is permissible given that 
(\ref{eq:desired_limit_range_1}) holds for any positive integer $R$.

\subsection{The case where $\left \lfloor \frac{\mu_r n}{2 \beta c \log n} \right \rfloor +1 \leq \ell \leq \min\{L_n,\lfloor \nu n \rfloor\}$}
\label{subsec:range3}
 We now consider the range where $\frac{\mu_r n}{2 \beta c \log n} < \ell \leq \min\{L_n,\lfloor \nu n \rfloor\}$
 for some $0<\nu<1/2$ to be specified later; see (\ref{eq:specifying_nu}).
From (\ref{eq:X_S_theta}), we see that
$X_{\ell,n} = \lfloor \beta \ell K_{1,n} \rfloor$ so that  (\ref{eq:bounding_expectation}) still holds.
Using (\ref{eq:bounding_expectation})
 and (\ref{eq:bound_K1Kr_P_detailed}) on the given range, we get
 \begin{align}
\bE{e^{-\frac{(X_{\ell,n}+1)|\Sigma|_n}{P_n}}}^{n-\ell} &\leq \left(1-\mu_r  
+ \mu_r e^{-\frac{\beta \ell K_{1,n}K_{r,n}}{P_n}}\right)^{n-\ell} 
\nonumber \\
& \leq \left(1-\mu_r  
+ \mu_r e^{-\beta \frac{\mu_r n}{2 \beta c \log n} \frac{c}{2} \frac{\log n}{n}} \right)^{n/2} 
\nonumber \\
& = \left(1-\mu_r  
+ \mu_r e^{- \frac{\mu_r}{4}  } \right)^{n/2} 
\label{eq:bounding_expectation_2}
\end{align}
for all $n$ sufficiently large. 
 Also, since $n \choose \ell$ is monotone increasing in $\ell$ over the range
$0 \leq \ell \leq \lfloor n/2\rfloor$ we have 
\begin{align}
{n \choose \ell} \leq {n \choose \lfloor \nu n \rfloor} \leq \left(\frac{e}{\nu}\right)^{\nu n}
\label{eq:range_3_inter_1}
\end{align}
by means of (\ref{eq:CombinatorialBound1}).

Using (\ref{eq:bounding_expectation_2}) and (\ref{eq:range_3_inter_1}) we now find
\begin{align}
& \sum_{\ell=\left \lfloor \frac{\mu_r n}{2 \beta c \log n} \right \rfloor +1}^{ \min\{L_n,\lfloor \nu n \rfloor\}}{n \choose \ell} \min\left\{1, \ell^{\ell-2} (p_{rr}(n))^{\ell-1}\right\}
\left(\min\left\{1-\lambda_1(n),\bE{e^{-\frac{(X_{\ell,n}+1)|\Sigma|_n}{P_n}}} \right\}\right)^{n-\ell}
\nonumber \\
& ~ \leq n  \left(\frac{e}{\nu}\right)^{\nu n} \left(1-\mu_r  
+ \mu_r e^{- \frac{\mu_r}{4}  } \right)^{n/2} 
\nonumber \\
& = n \left(\left(\frac{e}{\nu}\right)^{\nu}  \left(1-\mu_r  
+ \mu_r e^{- \frac{\mu_r}{4}  } \right)^{1/2}  \right)^n,
\label{eq:range_3_inter_2}
\end{align}
for all $n$ sufficiently large.
With $\mu_r>0$, it always holds that 
\[
\left(1-\mu_r  
+ \mu_r e^{- \frac{\mu_r}{4}  } \right)^{1/2} <1,
\]
while we have
\[
\lim_{\nu \to 0} \left(\frac{e}{\nu}\right)^{\nu}  = 1.
\]
Therefore, for any $\mu_r>0$, $\nu$ can be selected small enough to ensure that.
\begin{align}
\left(\left(\frac{e}{\nu}\right)^{\nu}  \left(1-\mu_r  
+ \mu_r e^{- \frac{\mu_r}{4}  } \right)^{1/2}  \right) < 1.
\label{eq:specifying_nu}
\end{align}
With $\nu$ selected according to (\ref{eq:specifying_nu}), we immediately obtain
 \begin{align}
\lim_{n \to \infty} \sum_{\ell=\left \lfloor \frac{\mu_r n}{2 \beta c \log n} \right \rfloor +1}^{ \min\{L_n,\lfloor \nu n \rfloor\}}{n \choose \ell} \min\left\{1, \ell^{\ell-2} (p_{rr}(n))^{\ell-1}\right\}
\left(\min\left\{1-\lambda_1(n),\bE{e^{-\frac{(X_{\ell,n}+1)|\Sigma|_n}{P_n}}} \right\}\right)^{n-\ell} = 0
\label{eq:desired_limit_range_3}
\end{align}
in view of (\ref{eq:range_3_inter_2}).

\subsection{The case where $\min\{L_n, \lfloor \nu n \rfloor\} + 1 \leq \ell \leq L_n$}
\label{subsec:range4}
Our next goal is to handle the range $\min\{L_n, \lfloor \nu n \rfloor\} + 1 \leq \ell \leq L_n$,
where it still holds that $X_{\ell,n} = \lfloor \beta \ell K_{1,n} \rfloor$. This range will become obsolete if
$L_n \leq \lfloor \nu n \rfloor$, so we only consider the case where $\lfloor \nu n \rfloor < L_n$, and hence the range
$\lfloor \nu n \rfloor +1 \leq \ell \leq L_n$.
On this range, we use the crude bound
$|\Sigma|_n \geq K_{1,n}$ to get
 \begin{align}
\bE{e^{-\frac{(X_{\ell,n}+1)|\Sigma|_n}{P_n}}}^{n-\ell} &\leq \left(e^{-\frac{\beta \ell K_{1,n}^2}{P_n}}\right)^{n-\ell} 
\leq e^{-\frac{\beta \ell K_{1,n}^2}{P_n}\frac{n}{2}} \leq  e^{-\frac{\beta \nu K_{1,n}^2}{2P_n}n^2},
\label{eq:bounding_expectation_3}
\end{align}
where we also note that $n-\ell \geq \frac{n}{2}$. Using (\ref{eq:bounding_expectation_3}) and (\ref{eq:CombinatorialBound2}) we now get
\begin{align}
& \sum_{\ell=\min\{L_n, \lfloor \nu n \rfloor\} + 1}^{L_n}{n \choose \ell} \min\left\{1, \ell^{\ell-2} (p_{rr}(n))^{\ell-1}\right\}
\left(\min\left\{1-\lambda_1(n),\bE{e^{-\frac{(X_{\ell,n}+1)|\Sigma|_n}{P_n}}} \right\}\right)^{n-\ell}
\nonumber \\
& ~ \leq \left( \sum_{\ell=\min\{L_n, \lfloor \nu n \rfloor\} + 1}^{L_n}{n \choose \ell} \right)   e^{-\frac{\beta \nu K_{1,n}^2}{2P_n}n^2}
\nonumber \\
& ~ \leq  \left( 2 e^{-\frac{\beta \nu K_{1,n}^2}{2P_n}n}\right)^{n}.
\label{eq:range_4_inter_1}
\end{align}
Under the enforced condition (\ref{eq:condition_for_con_2}) on the scaling we have that $\frac{K_{1,n}^2}{P_n}n = w(1)$, so that
\[
2 e^{-\frac{\beta \nu K_{1,n}^2}{2P_n}n} = o(1)
\]
since $\beta, \nu > 0$. Reporting this into (\ref{eq:range_4_inter_1}), we immediately obtain
\begin{align}
\lim_{n \to \infty} \sum_{\ell=\min\{L_n, \lfloor \nu n \rfloor\} + 1}^{L_n }{n \choose \ell} \min\left\{1, \ell^{\ell-2} (p_{rr}(n))^{\ell-1}\right\}
\left(\min\left\{1-\lambda_1(n),\bE{e^{-\frac{(X_{\ell,n}+1)|\Sigma|_n}{P_n}}} \right\}\right)^{n-\ell} = 0
\label{eq:desired_limit_range_4}
\end{align}

\subsection{The case where $L_n +1 \leq \ell \leq \left \lfloor \frac{n}{2} \right \rfloor$}
\label{subsec:range5}

Finally, we consider the range $L_n +1 \leq \ell \leq \left \lfloor \frac{n}{2} \right \rfloor$,
where we have $X_{\ell,n} = \lfloor \gamma P_n \rfloor$ as stated in (\ref{eq:X_S_theta}).
Using once again the crude bound
$|\Sigma|_n \geq K_{1,n}$ we get
 \begin{align}
\bE{e^{-\frac{(X_{\ell,n}+1)|\Sigma|_n}{P_n}}}^{n-\ell} &\leq \left(e^{-\gamma  K_{1,n}}\right)^{n-\ell} 
\leq e^{-\gamma K_{1,n}\frac{n}{2}}.
\label{eq:bounding_expectation_4}
\end{align}
In view of (\ref{eq:bounding_expectation_4}) and (\ref{eq:CombinatorialBound2}) we find
\begin{align}
& \sum_{\ell= L_n+1}^{\left \lfloor \frac{n}{2} \right \rfloor}{n \choose \ell} \min\left\{1, \ell^{\ell-2} (p_{rr}(n))^{\ell-1}\right\}
\left(\min\left\{1-\lambda_1(n),\bE{e^{-\frac{(X_{\ell,n}+1)|\Sigma|_n}{P_n}}} \right\}\right)^{n-\ell}
\nonumber \\
& ~ \leq \left( \sum_{\ell= L_n+1}^{\left \lfloor \frac{n}{2} \right \rfloor}{n \choose \ell} \right)   
e^{-\gamma K_{1,n}\frac{n}{2}}
\nonumber \\
& ~ \leq  \left( 2 e^{-\frac{\gamma K_{1,n}}{2}}\right)^{n}.
\label{eq:range_5_inter_1}
\end{align}
As before, we have under (\ref{eq:condition_for_con_1}) and
(\ref{eq:condition_for_con_2}) that $K_{1,n}=w(1)$ leading to 
\[
2 e^{-\frac{\gamma K_{1,n}}{2}} = o(1).
\]
Reporting this into (\ref{eq:range_5_inter_1}), we immediately obtain
\begin{align}
\lim_{n \to \infty} \sum_{\ell= L_n+1}^{ \left \lfloor \frac{n}{2} \right \rfloor}{n \choose \ell} \min\left\{1, \ell^{\ell-2} (p_{rr}(n))^{\ell-1}\right\}
\left(\min\left\{1-\lambda_1(n),\bE{e^{-\frac{(X_{\ell,n}+1)|\Sigma|_n}{P_n}}} \right\}\right)^{n-\ell} = 0.
\label{eq:desired_limit_range_5}
\end{align}

 \section*{Acknowledgements} 
 This work has been supported in part by the Department of Electrical and Computer Engineering
at Carnegie Mellon University. The author also thanks Prof. A. M. Makowski from UMD for insightful comments concerning this work.
 
 
 \section*{\LARGE Appendix}
\appendix
\setcounter{equation}{0}
\renewcommand{\theequation}{\thesection.\arabic{equation}}

\section{A proof of Corollary \ref{cor:node_isolation}}
 \label{sec:proof_cor}
In this section, we will show that combined together Theorem \ref{thm:node_isolation} 
and Theorem \ref{thm:connectivity}
is {\em equivalent}
to
Corollary \ref{cor:node_isolation} under the enforced assumptions.
Consider a probability distribution $\boldsymbol{\mu}=(\mu_1,\ldots,\mu_r)$ with $\mu_i>0$ for
all $i=1,\ldots, r$ and a scaling  $K_1,\ldots,K_r, P: \mathbb{N}_0 \to \mathbb{N}_0^{r+1}$.
The aforementioned equivalence of the results will follow upon showing the 
equivalence of the conditions  (\ref{eq:scaling_law}) and (\ref{eq:scaling_law_cor}), namely that for any $c>0$ 
we have
\[
\lambda_1(n) \sim c \frac{\log n}{n} \quad \textrm{if and only if} \quad \frac{K_{1,n} \bE{|\Sigma|_n}}{P_n} 
\sim c \frac{\log n}{n}.
\] 
In order to establish this, we will show that either of the conditions (\ref{eq:scaling_law}) or (\ref{eq:scaling_law_cor})
individually imply $ \lambda_1(n) \sim \frac{K_{1,n} \bE{|\Sigma|_n}}{P_n}$, or equivalently that
\begin{align}
\sum_{j=1}^{r} p_{1j}(n) \mu_j   &\sim 
\sum_{j=1}^{r} \frac{K_{1,n} K_{j,n}}{P_n}  \mu_j.
\label{eq:inter_appendix}
\end{align}
Since $\mu_j>0$ for all $j = 1,\ldots, r$, (\ref{eq:inter_appendix})
will follow immediately if we show under either (\ref{eq:scaling_law}) or (\ref{eq:scaling_law_cor}) that
\begin{equation}
p_{1j}(n) \sim \frac{K_{1,n} K_{j,n}}{P_n}, \qquad j=1,\ldots, r.
\label{eq:cor_inter}
\end{equation}
 Lemma \ref{lem:AsymptoticEquivalence} readily gives (\ref{eq:cor_inter}) as we note
 that for all $j=1,\ldots, r$, (\ref{eq:scaling_law}) implies $p_{1j}(n)=o(1)$ 
while  
 (\ref{eq:scaling_law_cor}) implies
 $\frac{K_{1,n} K_{j,n}}{P_n}=o(1)$. 
The equivalence of Theorem \ref{thm:node_isolation}-Theorem \ref{thm:connectivity}   
with Corollary \ref{cor:node_isolation} is now established.
\myendpf

\bibliographystyle{abbrv}
\bibliography{related}

\end{document}